\documentclass[letterpaper, 12pt, journal, onecolumn, draftclsnofoot]{IEEEtran}
\usepackage{cite}
\usepackage{amsmath,amssymb,amsfonts}
\usepackage{algorithmic}
\usepackage{graphicx}
\usepackage{textcomp}

\usepackage{url}
\usepackage{graphicx}
\usepackage{color}
\usepackage{cite}
\usepackage[hidelinks]{hyperref}
\usepackage{amsfonts}
\usepackage{amsmath,multirow} %
\usepackage{algorithm}
\usepackage{mathtools}

\usepackage{tikz}

\def\algname/{GoPRONTO}
\def\gradalgname/{Gradient \algname/}
\def\HBalgname/{Heavy-Ball \algname/}
\def\CGalgname/{Conjugate \algname/}
\def\nestalgname/{Nesterov's \algname/}
\def\projectionop/{PO}

\newcommand{\norm}[1]{\left \|#1 \right \|}
\newcommand{\R}{{\mathbb{R}}} % real
\newcommand{\N}{{\mathbb{N}}} % natural
\newcommand{\CC}{{\mathcal{C}}} % continuous functions space
\newcommand{\LL}{{\mathcal{L}}} % Lagrangian
\newcommand{\PP}{{\mathcal{P}}} % projection operator
\newcommand{\TT}{{\mathcal{T}}} % trajectory manifold
\newcommand{\tanTTk}{{T}_{(\bx,\bu)}^k{\TT}}
\newcommand{\tanTT}{{T}_{(\bx,\bu)}{\TT}}

\newcommand{\HB}{{\textsc{hb}}}

\newcommand\oprocendsymbol{\hbox{$\square$}}
\newcommand\oprocend{\relax\ifmmode\else\unskip\hfill\fi\oprocendsymbol}
 
\DeclareMathOperator*{\subj}{subj. to}

\newtheorem{theorem}{Theorem}[section]

\newtheorem{definition}[theorem]{Definition} 
\newtheorem{remark}[theorem]{Remark}

\newtheorem{assumption}[theorem]{Assumption}

\newcommand{\LS}[1]{\textcolor{black}{#1}}%

\DeclareMathOperator{\col}{col}

\newcommand{\horizon}[2]{[#1,#2]}
\newcommand{\map}[3]{#1: #2 \rightarrow #3}

\newcommand{\T}{^\top}

\newcommand{\dimx}{{n}}
\newcommand{\dimu}{{m}}
\newcommand{\dimalpha}{{\dimx T}}
\newcommand{\dimmu}{{\dimu T}}

\newcommand{\xinit}{x_{\text{init}}}

\newcommand{\bx}{\mathbf{x}}
\newcommand{\bu}{\mathbf{u}}

\newcommand{\bdeltax}{\boldsymbol{\Delta}\mathbf{x}}
\newcommand{\bdeltau}{\boldsymbol{\Delta}\mathbf{u}}
\newcommand{\deltax}{\Delta x}
\newcommand{\deltau}{\Delta u}
\newcommand{\dalpha}{\Delta \alpha}
\newcommand{\dmu}{\Delta \mu}
\newcommand{\tdalpha}{\Delta \tilde{\alpha}}
\newcommand{\tdmu}{\Delta \tilde{\mu}}

\newcommand{\blambda}{\boldsymbol{\lambda}}
\newcommand{\balpha}{\boldsymbol{\alpha}}
\newcommand{\bmu}{\boldsymbol{\mu}}

\newcommand{\tlambda}{\tilde{\lambda}}
\newcommand{\tbalpha}{\tilde{\balpha}}
\newcommand{\tbmu}{\tilde{\bmu}}

\let\min\relax
\DeclareMathOperator*{\min}{min}
\allowdisplaybreaks

\newcommand{\cost}{\ell}
\newcommand{\stagecost}{\ell_t}
\newcommand{\termcost}{\ell_T}
\newcommand{\dynamics}{f}

\newcommand{\impldyn}{{h}}

\newcommand{\impldynCL}{\tilde{h}} % Closed-loop

\newcommand{\stepsize}{\gamma}
\newcommand{\stepsizek}{\gamma^k}

\newcommand{\CGstep}{\rho}

\newcommand{\xreft}{x_{\text{ref},t}}
\newcommand{\ureft}{u_{\text{ref},t}}

\newcommand{\xmap}{\phi}
\newcommand{\umap}{\psi}

\newcommand{\friction}{{\mathrm{f}}}

\newcommand{\quadPRONTO}{W(\bx^k,\bu^k)}

\graphicspath{{figs/}}

\author{Lorenzo Sforni, Sara Spedicato, Ivano
  Notarnicola, Giuseppe
  Notarstefano \thanks{The authors are with the
    Department of Electrical, Electronic and Information Engineering, Alma Mater
    Studiorum - Universit\`{a} di Bologna, Bologna, Italy,
    \texttt{name.lastname@unibo.it}.}  \thanks{This result is part of a project
    that has received funding from the European Research Council (ERC) under the
    European Union's Horizon 2020 research and innovation programme (grant
    agreement No 638992 - OPT4SMART).  A very preliminary version of the idea
    inspiring this work, customized for the distributed framework, is proposed
    in~\cite{spedicato2018cloud}.}  }

\title{\algname/: a Feedback-based Framework\\for Nonlinear Optimal Control}

\begin{document}

\maketitle

\begin{abstract}
	% Nonlinear optimal control problems arise in a vast number of engineering applications 
	% in the automation and robotics fields.
	%especially when dealing with autonomous systems.
	%
	We propose GoPRONTO, a first-order, feedback-based approach to solve nonlinear discrete-time optimal control problems. 
	This method is a generalized first-order framework based on \LS{incorporating} 
	% a suitable embedding of 
	the original dynamics into a closed-loop system. 
	By exploiting this feedback-based shooting, we are able to reinterpret the optimal control problem as the minimization of a cost function, depending on a state-input curve, whose gradient can be computed by resorting to a suitable costate equation. This convenient reformulation gives room for a collection of accelerated numerical optimal control schemes.
	To corroborate the theoretical results, numerical simulations on the optimal control of  a train of inverted pendulum-on-cart systems are shown.
	% Taking inspiration from Hauser's PRojection Operator Newton method for Trajectory Optimization (PRONTO), we develop Go-PRONTO,
	% based on Conjugate gradient, Heavy-ball, and Nesterov's accelerated gradient. 
	%  An interesting original feature of GoPRONTO is that it does not require to solve quadratic optimization problems, so that it is well suited for the resolution of optimal control problems involving large-scale systems. 
	
\end{abstract}

\begin{IEEEkeywords}
Numerical optimal control, nonlinear control systems, trajectory optimization, gradient methods.
%send a blank e-mail to keywords@ieee.org or visit \underline
%{http://www.ieee.org/organizations/pubs/ani\_prod/keywrd98.txt}
\end{IEEEkeywords}

%%%%%%%%%%%%%%%%%%%%%%%%%%%%%%%%%%%%%%%%%%%%%%%%%%%%%%%%%%%%%%%%%%%%%%%%%%%%%%%
%%%%%%%%%%%%%%%%%%%%%%%%%%%%%%%%%%%%%%%%%%%%%%%%%%%%%%%%%%%%%%%%%%%%%%%%%%%%%%%

\section{Introduction}
\label{sec:intro}
\IEEEPARstart{N}{}onlinear optimal control problems arise in a vast number of engineering applications in Automation and Robotics.
%
% Many
% tasks, especially concerning autonomous systems, require the generation of a
% physically realizable trajectory, namely a trajectory that satisfies the system
% dynamics while optimizing a given performance criterion.
%
% Next we review some of the main numerical methods for nonlinear optimal control
% and highlight the main novelty of the framework proposed in this note.
% \sout{Next we review some numerical methods for nonlinear optimal control
% and highlight the novelty of the framework proposed in this note.}
% 
%%%%%%%%%%%%%%%%%%%%%%%%%%%%%%%%%%%%%%%%%%%%%%%%%%%%%%%%%%%%%%%%%%%%%%%%%%%%%%%
%%%%%%%%%%%%%%%%%%%%%%%%%%%%%%%%%%%%%%%%%%%%%%%%%%%%%%%%%%%%%%%%%%%%%%%%%%%%%%%
%
% \paragraph*{Literature Review} 
% Many approaches to solve optimal control problems have been
% investigated in the literature.
%
% \emph{Indirect methods} aim to satisfy the necessary conditions of optimality and
% typically solve a (two-point) boundary-value problem
% arising from calculus of variations (see, e.g., \cite{Sage68,Kirk70}) or
% from Pontryagin's Maximum Principle (see, e.g., \cite{Pontryagin62}).
% %
% Among the most recent works in this field we refer
% to~\cite{sassano2020combining,park2020optimal}.

In this note, we focus on \emph{direct methods} for the resolution of discrete-time optimal control problems.
These methods are subclassified into two different categories: \emph{simultaneous} and \emph{sequential} (see~\cite{diehl2009efficient}).
In simultaneous approaches, both the controls and the states are treated as decision variables of the nonlinear program (NLP), obtained via, e.g., collocation or multiple shooting methods from the original optimal control problem formulation.
The NLP is then addressed solving directly the Karush-Kuhn-Tucker optimality conditions of the problem by Newton's type optimization algorithms like Sequential Quadratic Programming (SQP)~\cite{diehl2009efficient} 
and Interior Point optimization (IP)~\cite{wright1997primal,Boyd04}.
These methods have been implemented in a variety of toolboxes, e.g.,~\cite{houska2011acado,englert2019software,verschueren2022acados, biegler2009large, zanelli2020forces}.
%
%
% \sout{The two major families of Newton type optimization methods are 
% Sequential Quadratic Programming (SQP) and Interior Point optimization (IP).
% %
% The literature on SQP is quite vast and we refer the interested reader to
% % \cite{Nocedal99,diehl2009efficient} 
% for detailed overviews.
% %
% SQP methods for the resolution of optimal control problems
% have been employed in various applications. For example, 
% % in~\cite{zanon2017asynchronous}, 
% an algorithm
% based on SQP is proposed to solve a vehicle coordination optimal control problem.
% %
% For IP methods, instead, we refer 
% % to~\cite{wright1997primal,Boyd04}.
% %
% Widely adopted implementations of nonlinear IP methods are represented by the 
% toolboxes 
% % IPOPT~\cite{biegler2009large} 
% and the more recent 
% FORCES
% % ~\cite{zanelli2020forces}.
% }
%
% open-source package IPOPT~\cite{biegler2009large} while a more recent toolbox is 
% FORCES~\cite{zanelli2020forces}.
%
While these approaches do not exhibit numerical instability issues coming from the integration of the dynamics, an important drawback is that, in general, the constraints describing the dynamics are satisfied only asymptotically.
\LS{That is, suboptimal state-input curves do
not satisfy the dynamics in general.}
% they do not enjoy a ``dynamic
% feasibility''. 
% That is, 
% the state-input curves computed at each iteration satisfies 
%
% \sout{That is, the state-input curves computed at each iteration do
% not satisfy the dynamics in general.}
%
% However, this feature can be extremely
% important in real-time control schemes 
% % (as, e.g., in Model Predictive Control)
% since it may allow for suboptimal schemes stopping after few iterations.
% Approaches to deal with feasibility of the dynamic constraints in SQP methods
% have been presented in~\cite{Tenny04} and in~\cite{bayer2013projected}.
%
Conversely, sequential approaches tackle the NLP in the
reduced space of control variables only.  
Hence, at each iteration, 
the state trajectory
is recovered by forward simulation of the system
dynamics.
A first-order sequential approach for the resolution of nonlinear optimal
control problems is presented, e.g.,
in~\cite[Section~1.9]{bertsekas1999nonlinear}.
The main limitation of sequential methods is their numerical instability in the forward simulation of (the possibly unstable) dynamics.
%
% These issues, indeed, lead to the development of multiple-shooting methods (see, e.g.,~\cite{bock1984multiple}) to replace single shooting methods (a sequential approach).
%
This issue is overcome by the original sequential method \emph{PRojection Operator Newton
Method for Trajectory Optimization (PRONTO)}~\cite{hauser2002projection}.
In this work, through the use of a control feedback, the shooting map is stabilized and the
optimal control problem is converted into an unconstrained one to which a Newton's
method is applied.
This approach has been successfully extended in a variety of scenarios~\cite{hauser2006barrier, saccon2013optimal, aguiar2017constrained, filo2018function, bayer2013projected}.
The introduction of stabilizing controllers to handle unstable dynamics is considered also in the model predictive control literature, see, e.g.,~\cite{kouvaritakis2000efficient,aftab2021pre}.
Finally, optimal control problems have been also tackled via augmented Lagrangian approaches in~\cite{kouzoupis2016block}.

%%%%%%%%%%%%%%%%%%%%%%%%%%%%%%%%%%%%%%%%%%%%%%%%%%%%%%%%%%%%%%%%%%%%%%%%%%%%%%%
%%%%%%%%%%%%%%%%%%%%%%%%%%%%%%%%%%%%%%%%%%%%%%%%%%%%%%%%%%%%%%%%%%%%%%%%%%%%%%%

% \paragraph*{Contributions}

The main contribution of this note is to provide a novel 
class of numerically-robust first-order algorithms for \LS{(offline)} discrete-time optimal control 
termed  \algname/, short for Generalized first-Order
PROjectioN operator method for Trajectory Optimization\footnote{This acronym is
  chosen as a tribute to Hauser's PRONTO.}.
In our approach we combine the introduction of a feedback
system (projection operator) into the nonlinear optimal control problem 
with a gradient-based resolution strategy.
Such innovative combination results into an optimization framework that enjoys several
appealing features.
Thanks to the introduction of the projection operator in the sequential optimization scheme, we achieve (i) numerical robustness, even when dealing with unstable dynamics, and (ii) dynamic feasibility, i.e., 
% a state-input trajectory is available at each iteration.
\LS{a state-input trajectory can be computed, at each iteration, via a closed-loop integration of the nonlinear dynamics. This property is of particular interest in the case of, e.g., unstable systems.}
From the (open-loop) gradient method for optimal control in
\cite{bertsekas1999nonlinear}, our approach inherits the simplicity of
implementation of the descent-direction search, namely a costate equation update.  
Finally, we show how \algname/ gives rise to several first-order
optimization algorithms that can speed up the resolution of the optimal control problem.
This simple and adaptable update rule makes \algname/ flexible enough to be extended to problems
involving large-scale dynamics.
As in other optimization domains with very-large decision variables, Newton's methods are impracticable while first-order
approaches are preferred.

% and its
% variants~\cite{kingma2014adam, reddi2018convergence} for further details.
% the policy gradient method in the field of Reinforcement Learning,
% for further details see~\cite{sutton2018reinforcement}.
%
%ADAM?
%
%
% Indeed, novel algorithms can be implemented by means of appropriate
% modifications in the structure of our algorithmic approach.
% %
% Finally, we show how this general \LS{framework} gives rise to several first-order
% optimization algorithms that can speed up the resolution of the optimal control
% problem.
% %
% Indeed, novel algorithms can be implemented by means of appropriate
% modifications in the structure of our algorithmic approach.
% , e.g., by variations
% in the state-input curve update or by the evaluation of the costate equation and
% the system-linearization in different curves.
%
% Here, we consider three alternative first-order optimization algorithms, namely the
% conjugate gradient~\cite{hestenes1952methods}, Heavy-Ball \cite{polyak1964some}
% and Nesterov's Accelerated Gradient \cite{nesterov1983method}, and propose their
% counterparts
% % i.e., \CGalgname/, \HBalgname/ and \nestalgname/ 
% in our optimal control framework.
%

%%%%%%%%%%%%%%%%%%%%%%%%%%%%%%%%%%%%%%%%%%%%%%%%%%%%%%%%%%%%%%%%%%%%%%%%%%%%%%%
%%%%%%%%%%%%%%%%%%%%%%%%%%%%%%%%%%%%%%%%%%%%%%%%%%%%%%%%%%%%%%%%%%%%%%%%%%%%%%%

% \paragraph*{Organization}

The note unfolds as follows. 
The nonlinear optimal control problem is presented in Section~\ref{sec:setup} along with some preliminaries.
% existing numerical methods for optimal control.
% as the gradient method for optimal control presented in \cite{bertsekas1999nonlinear} 
% and PRONTO~\cite{hauser2002projection}. 
%
In Section~\ref{sec:emb_feedback_method} we propose our methodology \algname/
in its steepest descent implementation.
Accelerated versions of \algname/ are presented in Section~\ref{sec:accelerated}.
Some numerical simulations with a train of inverted pendulum-on-cart systems are given in
Section~\ref{sec:simulations}.

% followed by some concluding remarks in Section \ref{sec:conclusions}.

%%%%%%%%%%%%%%%%%%%%%%%%%%%%%%%%%%%%%%%%%%%%%%%%%%%%%%%%%%%%%%%%%%%%%%%%%%%%%%%
%%%%%%%%%%%%%%%%%%%%%%%%%%%%%%%%%%%%%%%%%%%%%%%%%%%%%%%%%%%%%%%%%%%%%%%%%%%%%%%

\paragraph*{Notation}
% Throughout the note, all the matrices and vector are assumed to have compatible dimensions.
%
The vertical stack of $x_1$ and $x_2$ is denoted by $\col(x_1,x_2) := [x_1^\top, x_2^\top]^\top$.
%
% Given a (scalar) function $\map{\ell}{\R^n}{\R}$, the gradient of $\ell$ is the
% $n \times 1$ column vector defined as
% %
% $\nabla \ell(x) := [\tfrac{\partial \ell (x)}{\partial x_1}\T, \ldots, \tfrac{\partial \ell (x)}{\partial x_n}\T]\T$. 
%
% \begin{align*}
%   \nabla \ell(x) =
%   \begin{bmatrix}
%     \frac{\partial \ell (x)}{\partial x_1}
%     \\
%     \vdots
%     \\
%     \frac{\partial \ell (x)}{\partial x_n}
%   \end{bmatrix}
%   \in \R^{n \times 1}
% \end{align*}
%
Given $\map{\ell}{\R^{\dimx} \times \R^{\dimu}}{\R}$, its (total) gradient at a given point $(\bar{x}, \bar{u})$ is
$\nabla \ell( \bar{x}, \bar{u} ) := \col(\nabla_{x} \ell( \bar{x},
\bar{u}),\nabla_{u} \ell( \bar{x}, \bar{u}) )$, where
$\nabla_{x} \ell( \bar{x}, \bar{u} )$ 
%:= \nabla_{x}
%  \ell(x,u)\big|_{x=\bar{x}, u = \bar{u}}$,
% \begin{align*}
% 	\nabla_{x} \ell( \bar{x}, \bar{u} )
%   & :=
%   \nabla_{x} \ell(x,u)\Big|_{x=\bar{x}, u = \bar{u}},
% \end{align*}
% i.e., the partial derivative with respect to the first argument of $\ell$,
and $\nabla_{u} \ell( \bar{x}, \bar{u} )$
are the partial derivative of $\ell$ with respect to the first and the second argument, respectively.
%to the second argument.
%
Moreover, given $\map{f}{\R^n}{\R^m}$, the gradient of $f$ is
$ \nabla f (x) := [\nabla f_1 (x) \ldots \nabla f_m (x)] \in \R^{n \times m} $.
For $T \in \N$, we define $\horizon{0}{T} := \{0,1,2,\dots,T\}$.
\LS{Given a symmetric, positive-definite matrix $Q\in \R^{\dimx \times \dimx}$, and $x \in \R^\dimx$, we define the $Q$-norm of $x$ as $\| x \|_Q = \sqrt[2]{x\T Q x}$}
%

%%%%%%%%%%%%%%%%%%%%%%%%%%%%%%%%%%%%%%%%%%%%%%%%%%%%%%%%%%%%%%%%%%%%%%%%%%%%%%%
%%%%%%%%%%%%%%%%%%%%%%%%%%%%%%%%%%%%%%%%%%%%%%%%%%%%%%%%%%%%%%%%%%%%%%%%%%%%%%%

\section{Problem Setup and Preliminaries}
\label{sec:setup}

In this section the nonlinear, discrete-time optimal control setup investigated
in the note is introduced.
Then, we review two numerical methods for optimal control related to the
proposed approach, 
% namely the gradient method for optimal control
discussed in~\cite{bertsekas1999nonlinear} and 
in~\cite{hauser2002projection}.
% a discrete-time version of PRONTO, originally proposed 

\subsection{Discrete-time Optimal Control Setup}

We consider nonlinear, discrete-time systems
described by
% dynamics
%
\begin{align}\label{eq:dynamics}
	x_{t+1} &= \dynamics (x_{t}, u_{t})  \qquad t \in \N
\end{align}
where
$x_t \in \R^{\dimx}$ and $u_t \in \R^{\dimu}$ are the state and the input 
of the system at time $t$, respectively.
The map $\map{\dynamics}{\R^{\dimx} \times \R^{\dimu}}{\R^{\dimx}}$ is the vector
field describing the nonlinear dynamics.
The initial condition of the system is a fixed value $\xinit \in \R^{\dimx}$. %  its (given) inital condition,
%
% \begin{remark}
%   We point out that although ours is a discrete-time framework,
%   continuous-time systems can be also considered.
%   % 
%   As customary in the literature, continuous-time systems can be
%   discretized by means of proper integration schemes. 
% 	% As we will point out later
%   % in the note, the numerical robustness of the proposed approach allows one to
%   % easily use commonly available discretization schemes, since integration of
%   % closed-loop systems is required in the process. \oprocend
% \end{remark}
%

For notational convenience, we use $\bx \in \R^{\dimx T} $ and
$\bu \in \R^{\dimu T} $ to denote, respectively, the stack of the states $x_t$
for all $t \in \horizon{1}{T}$ and the inputs $u_t$ for all
$t \in \horizon{0}{T-1}$, that is $\bx := \col(x_1, \dots, x_T)$ and
$\bu := \col(u_0, \dots, u_{T-1})$.
% \begin{align*}
%   \bx &:= \col(x_1, \dots, x_T)
%   \\
%   \bu &:= \col(u_0, \dots, u_{T-1})
% \end{align*}
%
% i.e., defined as
%
Next we give a useful definition.

% On an horizon $T$ we define as a trajectory
%
\begin{definition}[Trajectory]\label{def:trajectory}
	A pair $(\bx,\bu) \in \R^{\dimalpha} \times \R^{\dimmu}$ is called a \emph{trajectory} 
	of the system described by~\eqref{eq:dynamics} if its components satisfy the 
	constraint represented by the dynamics~\eqref{eq:dynamics} for all $t\in\horizon{0}{T-1}$.
  In particular, $\bx$ is the state trajectory, while $\bu$ is the input trajectory.
  \oprocend
\end{definition}

Conversely, we refer to a generic pair $(\balpha,\bmu) \in \R^{\dimalpha} \times \R^{\dimmu}$ 
with \mbox{$\balpha := \col(\alpha_1, \ldots, \alpha_T)$} and $\bmu := \col(\mu_0,\ldots, \mu_{T-1})$
as a state-input \emph{curve},
in analogy with the continuous-time terminology.
Notice that a curve $(\balpha,\bmu)$ is not necessarily a trajectory, i.e., it
does not necessarily satisfy the dynamics~\eqref{eq:dynamics}.
% Clearly, $\balpha \in \R^{\dimalpha}$ can be thought to be the stack of $\alpha_t \in \R^{\dimx}$ for
% all $t \in \horizon{1}{T}$ and $\bmu \in \R^{\dimmu}$ the stack of $\mu_t \in \R^{\dimu}$ for all
% $t \in \horizon{0}{T-1}$, i.e.,
% \begin{align*}
	% 	\balpha &:= \col(\alpha_1, \ldots, \alpha_T)
	% 	\\
	% 	\bmu &:= \col(\mu_0,\ldots, \mu_{T-1}).
	% \end{align*}

% We denote by $\TT\subset \R^\dimalpha \times \R^\dimmu$ the trajectory manifold
% of the control system~\eqref{eq:dynamics}.
%
By rewriting the nonlinear dynamics~\eqref{eq:dynamics}
as an implicit equality constraint $\map{\impldyn}{\R^{\dimalpha} \!\times \!\R^{\dimmu}}{\R^{\dimalpha}}$
% given by
% %
% \begin{align}
%   \label{eq:compact_dynamics}
% 	\impldyn(\bx, \bu) :=
% 	\begin{bmatrix}
% 		\dynamics(x_0, u_0) - x_1 \\
% 		\vdots	\\
% 		\dynamics(x_{T-1},u_{T-1}) - x_T
% 	\end{bmatrix},
% \end{align}
%
we define the trajectory manifold $\TT\subset \R^\dimalpha \times \R^\dimmu$
of~\eqref{eq:dynamics} as the set $\TT := \{(\bx, \bu) \mid \impldyn(\bx,\bu) = 0\}$.
%
% \begin{align} \label{eq:traj_man}
% 	\TT := \{(\bx, \bu) \mid \impldyn(\bx,\bu) = 0\}.
% \end{align}
%
It can be shown that the tangent space to the trajectory
manifold
% ~\eqref{eq:traj_man} 
at a given trajectory (point), denoted as $\tanTT$,
is represented by the set of trajectories satisfying the linearization of the
nonlinear dynamics~$\dynamics(\cdot,\cdot)$ about the trajectory $(\bx,\bu)$.
	
Among all possible trajectories of system~\eqref{eq:dynamics}, we aim to
optimize a given performance criterion defined over a fixed time horizon
$\horizon{0}{T}$.  Formally, we look for a solution of the discrete-time optimal
control problem
\begin{subequations}
\label{eq:ocp_original}
\begin{align}\label{eq:ocp_original:cost}
  \min_{\bx \in \R^\dimalpha, \bu \in \R^\dimmu}
  \:&\:
  \sum_{t=0}^{T-1}
  \stagecost (x_{t},u_{t})
  +
  \termcost (x_{T})
  \\
  \label{eq:ocp_original:dynamics}
  \begin{split}
  \subj \: & \: x_{t+1} = \dynamics (x_{t}, u_{t}),   \quad t \in \horizon{0}{T-1}
  %\\
  %& \: x_{0} = \xinit,
	\end{split}
	% \\
	% & \review{\: \constr_{t}(x_t, u_t) \le 0}
	% \label{eq:ocp_original:constraints}
\end{align}
\end{subequations}
with initial condition $x_0 = \xinit \in \R^\dimx$, stage cost
$\map{\stagecost}{\R^\dimx \times \R^\dimu}{\R}$ and terminal cost
$\map{\termcost}{\R^\dimx}{\R}$.
%
% \review{The state-input constraints are represented by $\map{\constr_{t}}{\R^\dimx \times \R^{\dimu}}{\R^\nconstr}$.}
% % and $\map{f}{\R^\dimx \times \R^\dimu}{\R^\dimx}$. 

\begin{remark}
	% For the sake of exposition, we present the algorithm tailored for \emph{unconstrained} optimal control problems and neglect the presence of state-input constraints.
	%
\LS{For the sake of exposition, the algorithm presented in the paper is tailored for \emph{unconstrained} optimal control problems, neglecting the presence of state-input constraints.
	However, constraints can be addressed by adopting barrier function approaches (e.g.,\cite{hauser2006barrier}). Although these approaches may influence the numerical properties and performance of the algorithm, they have demonstrated successful results for Projection Operator-based algorithms in various settings (cf. the review paper\cite{aguiar2017constrained} and references therein).}
	\oprocend
\end{remark}

\begin{assumption}
  All functions $\stagecost(\cdot,\cdot)$, $\termcost(\cdot)$ and
  $\dynamics(\cdot,\cdot)$ are twice continuously differentiable, i.e., they are
  of class $\CC^2$.
	\oprocend
	\label{asm:regularity}
\end{assumption}

% By compactly defining the cost function~\eqref{eq:ocp_original:cost} as
% %
%   \begin{align}
%     % \label{eq:compact_cost}
%     \label{eq:cost_def}
% 	\cost(\bx, \bu) := \sum_{t=0}^{T-1}
% 	\stagecost (x_{t},u_{t})
% 	+
% 	\termcost (x_{T}),
% \end{align}
% %
% problem~\eqref{eq:ocp_original} can be rewritten as
% %
% \begin{align*}
% 	\begin{aligned}
% 		\min_{\bx \in \R^\dimalpha, \bu \in \R^\dimmu}
% 		\;&\; 
% 		\cost (\bx, \bu) 
% 		\\
% 		\subj \;&\; \impldyn (\bx, \bu) = 0
% 	\end{aligned}
% 	& \qquad\text{or} \quad &
% 	\begin{aligned}
% 		\min_{(\bx, \bu) \in \TT} % _{\bx \in \R^\dimalpha, \bu \in \R^\dimmu} 
% 		\;&\; 
% 		\cost (\bx, \bu) 
% 		% \\
% 		% \subj \;&\; (\bx, \bu) \in \TT
% 	\end{aligned}
% \end{align*}
%
% or, equivalently,
% %
% \begin{align*}
% 	\min_{\IN{(\bx, \bu) \in \R^\dimalpha\times\R^\dimmu}} 
% 	\;&\; 
% 	\cost (\bx, \bu) 
% 	\\
% 	\subj \;&\; (\bx, \bu) \in \TT
% \end{align*}
%

% In the following discussion, we are going to temporarily consider an \emph{unconstrained} optimal control framework, i.e., neglect the state-input constraints. 
% %
% The resulting optimal control problem is
% %
% \begin{align*}
% 	\begin{aligned}
% 		\min_{\bx \in \R^\dimalpha, \bu \in \R^\dimmu}
% 		\;&\; 
% 		\cost (\bx, \bu) 
% 		\\
% 		\subj \;&\; \impldyn (\bx, \bu) = 0
% 	\end{aligned}
% 	& \qquad\text{or} \quad &
% 	\begin{aligned}
% 		\min_{(\bx, \bu) \in \TT} % _{\bx \in \R^\dimalpha, \bu \in \R^\dimmu} 
% 		\;&\; 
% 		\cost (\bx, \bu) 
% 	\end{aligned}
% \end{align*}
% %

% It is worth noting 
Notice that, in light of the nonlinear equality constraint
% $\impldyn (\bx, \bu) = 0$ 
of the nonlinear dynamics,
problem~\eqref{eq:ocp_original} is a nonconvex program.

% Figure~\ref{fig:trajman_clean} provides a graphical representation of the optimal control problem
% as a nonlinear (nonconvex) program.
%
% \begin{figure}[htbp]
% 	\centering
% 	\includegraphics{trajman_clean.pdf}
% 	\caption{Two dimensional representation of the optimal control problem: 
% 	%
%             in gray the level curves of the cost function $\ell(\cdot,\cdot)$, in
%             black the nonlinear constraint representing the trajectory manifold
%             $\TT$, in green its tangent space $T_{(\bar{\bx}, \bar{\bu})}\TT$ about
%             $(\bar{\bx}, \bar{\bu})$.}
% 	\label{fig:trajman_clean}
% \end{figure}

Throughout the note, we will use the following shorthand notation for the linearization 
of both the cost function and dynamics about
a generic trajectory $(\bx^k, \bu^k)$ at iteration $k > 0$
\vspace{-2ex}
\begin{subequations}
	\label{eq:shorthand}
	\begin{align}
		a_t^k &:= \nabla_{x_t} \stagecost (x_t^k, u_t^k), \hspace{1.1cm}
		b_t^k := \nabla_{u_t} \stagecost (x_t^k, u_t^k),
		\label{eq:shorthand_ab}
		\\
		A_t^k &:= \nabla_{x_t} \dynamics(x_t^k, u_t^k)^\top, \hspace{0.8cm}
		B_t^k := \nabla_{u_t} \dynamics (x_t^k, u_t^k)^\top.
	\label{eq:shorthand_AB}
	\end{align}
\end{subequations}

%%%%%%%%%%%%%%%%%%%%%%%%%%%%%%%%%%%%%%%%%%%%%%%%%%%%%%%%%%%%%%%%%%%%%%%%%%%%%%%
%%%%%%%%%%%%%%%%%%%%%%%%%%%%%%%%%%%%%%%%%%%%%%%%%%%%%%%%%%%%%%%%%%%%%%%%%%%%%%%

\subsection{Gradient Method for Optimal Control}
\label{sec:open_loop_sequential}
Next, we recall a numerical strategy proposed, e.g., in~\cite[Section~1.9]{bertsekas1999nonlinear} 
to solve {a discrete-time optimal control problem as in~\eqref{eq:ocp_original} based on the gradient method.}

The leading idea is % to exploit the knowledge of the dynamics
to express the state $x_t$ at each $t \in \horizon{0}{T-1}$ as a function of $\bu$ only.
Formally, for all $t$ we can introduce a map\footnote{Formally, the shooting map depends also on the initial condition $\xinit$. Being it fixed, for notational convenience we drop this dependence.}
$\map{\xmap_t}{\R^{\dimmu}}{\R^{\dimx}}$ such that $x_t := \xmap_t (\bu)$,
%
% \begin{align}\label{eq:ol_map}
% 	x_t := \xmap_t (\bu),
% \end{align}
%
so that problem~\eqref{eq:ocp_original} can be recast into the % so-called
reduced version
\begin{align}
	\label{eq:unconstrained_problem}
%	\begin{split}
		% \min_{\bu \in \R^{\dimmu}} 
\min_{\bu} %\:
	% &
	%\:
	% \underbracket{
	\sum_{t=0}^{T-1}
	\stagecost (\xmap_{t}(\bu),u_{t})
	+
	\termcost (\xmap_{T}(\bu))
	% }
	=
	\min_{\bu} %\:
	J(\bu)
\end{align}
where the optimization variable is only the input sequence $\bu\in \R^{\dimmu}$.
Problem~\eqref{eq:unconstrained_problem} is an unconstrained 
optimization problem in $\bu$ 
%with a sufficiently regular cost function.
with a $\CC^2$ cost function.
Notice that the cost function $J(\bu)$
% in problem~\eqref{eq:unconstrained_problem}
inherits from~\eqref{eq:ocp_original} its smoothness properties, but also its
nonconvexity.
Hence, problem~\eqref{eq:unconstrained_problem} can be addressed via a gradient
descent method in which each component $u_t^k$, $t \in \horizon{0}{T-1}$ of $\bu^k \in \R^{\dimmu}$ is
iteratively updated as
\begin{align}\label{eq:barv}
	u_t^{k+1} & = u_t^k + \stepsizek 
	\deltau_t^k.
	% \underbrace{ 
	% \nabla_{u_t} J (\bu^k)
	% }_{\displaystyle -{\deltau}_t^k}
\end{align}
where $k > 0$ denotes the iteration counter, $\deltau_t^k := \nabla_{u_t} J (\bu^k)$, while
the parameter $\stepsizek >0$ is the so-called step-size.
%

%

% \begin{align}
% 	\bu^{k+1} = \bu^k - \stepsizek \nabla J (\bu^k),
% \end{align}
% %
% where $k > 0$ denotes the iteration counter, while
% the parameter $\stepsizek >0$ is the so-called step-size.

% Denoting $\deltau_t^k = - \nabla_{u_t} J (\bu^k)$, the previous update can be also written in a component-wise fashion 
% for $t \in \horizon{0}{T-1}$ as
% %
% \begin{align}\label{eq:barv}
% 	u_t^{k+1} & = u_t^k + \stepsizek 
% 	\deltau_t^k.
% 	% \underbrace{ 
% 	% \nabla_{u_t} J (\bu^k)
% 	% }_{\displaystyle -{\deltau}_t^k}
% \end{align}
%
%

The gradient of $J (\cdot)$ at every $\bu^k$ can be efficiently computed by
properly exploiting a costate difference equation (to be simulated backward in
time) based on the linearization of the cost and the system dynamics at a given
trajectory $(\bx^k, \bu^k)$ according to~\eqref{eq:shorthand}.
This backward pass reads, for each $k > 0$,
\begin{subequations}
	\label{alg:SOL:descent_direction}
		\begin{align}
			\label{alg:SOL:descent_direction_lambda}
			\lambda_t^k & = A_t^{k \top} {\lambda}_{t+1}^k + a_t^k
			\\
			\label{alg:SOL:descent_direction_v}
			\deltau_t^k & = - B_t^{k \top} {\lambda}_{t+1}^k - b_t^k.
		\end{align}
	\end{subequations}
As mentioned above, the update of the costate
$\blambda^k = \col ( {\lambda}_1^k, \ldots, {\lambda}_T^k)$ involves the
linearization of both the cost and the dynamics at the current input estimate
$\bu^k$ and corresponding state $\bx^k$ {(cf.~\eqref{eq:shorthand})}.
Then, the component ${\deltau}_t^k \in \R^\dimu$ of the 
{update (descent)} direction in~\eqref{eq:barv} is obtained via~\eqref{alg:SOL:descent_direction}.
%
%
% Thus, 
% \IN{The updated state trajectory $\bx^{k+1}$ is then computed by forward
% simulation of the dynamics~\eqref{eq:ocp_original:dynamics} over the horizon
% $\horizon{0}{T-1}$ with $x_{0}^{k+1} = \xinit$ with input $\bu^{k+1}$.}
The algorithm makes explicit use of the state sequence $\bx^{k}$ (associated to
the current input estimate $\bu^{k}$), which is obtained by forward
simulation of the dynamics~\eqref{eq:ocp_original:dynamics} over the horizon
$\horizon{0}{T-1}$.
, i.e., via 
\begin{align}
	\begin{split}
		% u_t^{k+1}&=u_t^k + \stepsizek\, \deltau_t^k
		% \\
		x_{t+1}^{k+1} & = f(x_t^{k+1}, u_t^{k+1})
	\end{split}
	\label{alg:SOL:trajectory_update}
\end{align}
with $x_{0}^{k+1} = \xinit$, so that $(\bx^{k+1}, \bu^{k+1})$ is a trajectory.
% {(cf.~\eqref{alg:SOL:trajectory_update})}.

\begin{remark}
  We stress that
	, as follows from~\eqref{alg:SOL:trajectory_update}, 
	each
  state trajectory $\bx^{k+1}$ is generated by an open-loop simulation of the
  dynamics, so that the method is not practically implementable for
  systems exhibiting unstable behaviors.\oprocend
\end{remark}
% \begin{remark}
%   We stress that
% 	% , as follows from~\eqref{alg:SOL:trajectory_update}, 
% 	each
%   state trajectory $\bx^{k+1}$ is generated by an open-loop simulation of the
%   dynamics, so that the method is not practically implementable for
%   systems exhibiting unstable behaviors.\oprocend
% \end{remark}

% Since Algorithm~\ref{alg:SOL} 
% \sout{This algorithm generates a sequence of inputs
% $\{ \bu^k\}_{k\ge 0}$ associated to a gradient method applied
% to~\eqref{eq:unconstrained_problem}, it inherits its convergence results.}
%
% Notice that the presented backward-forward sweep algorithm could be seen as 
% the reverse mode of the so-called Algorithmic Differentiation, see~\cite{griewank2008evaluating} 
% for details.
% 

%%%%%%%%%%%%%%%%%%%%%%%%%%%%%%%%%%%%%%%%%%%%%%%%%%%%

\subsection{Discrete-time PRONTO} \label{sec:pronto}

Along the lines of \cite{bayer2013projected}, we present a discrete-time version of the continuous-time optimal control algorithm PRONTO~\cite{hauser2002projection}.
%  in continuous-time.
%
% The key idea of PRONTO is to use a state feedback in an optimal control method
% to gain numerical stability and interpret this feedback as a projection operator
% that maps (state-input) curves into system trajectories.

The key idea of PRONTO is to use a stabilizing feedback in an optimal control method
to gain numerical stability, and to interpret such (tracking) controller as a projection operator
that maps (state-input) curves into system trajectories.
%
% By means of such projection operator, at each iteration a new (infeasible) state-input 
% curve is mapped to a (feasible) trajectory satisfying the dynamics. 
%
Given a state-input curve $(\balpha,\bmu)$, let us formally consider a nonlinear tracking system given by
% , for $t \in \horizon{0}{T-1}$,
\begin{align}
	\label{PR:closed_loop_projection}
	% \begin{split}
	  u_{t} & = \mu_{t} + K_{t} (\alpha_{t} - x_{t}),
	  &
	  x_{t+1} & = \dynamics (x_{t}, u_{t}),
	% \end{split}
 	% \qquad t \in \horizon{0}{T-1}
\end{align}
where $K_t \in \R^{\dimx \times \dimu}$ is a properly selected gain matrix.
\begin{remark}
	The feedback gain $K_t$ in~\eqref{PR:closed_loop_projection} should
	ensure local stability about the current state-input trajectory. 
	Among different alternatives, a \LS{possible} choice could be a linear quadratic regulator
	involving the system linearization about the current trajectory. 
	Also, one can use more advanced design approaches, e.g., linear parameter varying controllers~\cite{spedicato2019sparse}.%
	\oprocend
\end{remark}

System~\eqref{PR:closed_loop_projection} defines a nonlinear map, denoted by
$\map{\PP}{ \R^{\dimalpha} \times \R^{\dimmu}}{\TT}$ (with $\TT$ the trajectory
manifold)
%  defined in~\eqref{eq:traj_man}), 
 such that
\begin{align}\label{PR:projection_operator}
  \begin{bmatrix}
	  \balpha\\
	  \bmu
  \end{bmatrix} \longmapsto \begin{bmatrix}
	  \bx\\
	  \bu
  \end{bmatrix} := \PP ( \balpha,\bmu) = \begin{bmatrix}
	  \xmap(\balpha,\bmu) \\ \umap(\balpha,\bmu)
  \end{bmatrix},
\end{align}
where $\xmap(\balpha,\bmu)$ and $\umap(\balpha,\bmu)$ are the state and input
components of $\PP(\balpha,\bmu)$.
$\PP$ is a \emph{projection} since $(\bx,\bu) = \PP (\bx,\bu)$.
% state-input curves $(\balpha,\bmu)$ into trajectories (i.e., curves satisfying the dynamics) $(\bx,\bu)$, i.e., 

%

%---------- OLD
% It can be seen that this feedback system defines a nonlinear map
% $\map{\PP}{ \R^{\dimalpha} \times \R^{\dimmu}}{\R^{\dimalpha} \times
%   \R^{\dimmu}}$
% \begin{align*}
%   ( \balpha,\bmu) \longmapsto (\bx,\bu) = \PP ( \balpha,\bmu) \GN{=\Big(\phi(\alpha,\mu),\psi(\alpha,\mu)\Big)} 
% \end{align*}
% which is a projection mapping state-input curves $(\balpha,\bmu)$ into
% trajectories (i.e., curves satisfying the dynamics) $(\bx,\bu)$.

% \GN{Add here $\ell(\phi(\alpha,\mu),\psi(\alpha,\mu))$ with
%   $(\bdeltax,\bdeltau)$ trajectory of the linearization being tangent space to
%   space of trajectories.
% }

% Then, the strategy proposed in \cite{hauser2002projection} tackles this unconstrained problem 
% by means of a (projected) Newton method.
%

\begin{figure}[htbp]
	\centering
	\includegraphics{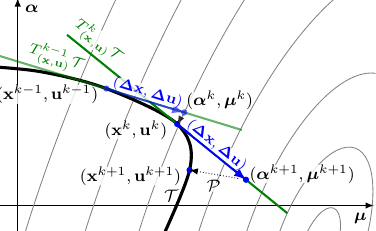}
	\caption{\normalfont Representation of PRONTO approach:  %
          in gray the level curves of the cost function $\cost(\cdot,\cdot)$, in
          black the trajectory manifold $\TT$, in green the tangent space $T^k_{(\bx, \bu)}\TT$ about trajectory $(\bx^k,\bu^k)$.  % $(\bx,\bu)$
          % The update direction is represented in blue.
          %
          At each iteration $k$, the update direction $(\bdeltax,\bdeltau)$ in
          blue is sought on the tangent space % to the trajectory manifold % $\TT$
          at the current trajectory $(\bx^k,\bu^k$). The updated curve
          $(\balpha^{k+1},\bmu^{k+1})$ is then projected onto $\TT$ by the
          projection operator $\PP$ (dotted line).
	%
	% For notational convenience, the update direction is denoted, for all $k$, as $(\deltax,\deltau)$ instead of $(\deltax^k,\deltau^k)$.  
	} 
	\label{fig:trajman_pronto}
\end{figure}

Thanks to the projection operator, the optimal control
  problem~\eqref{eq:ocp_original} can be written as
\begin{align}\label{PR:compact_problem}
	\min_{\balpha, \bmu} \:
	\cost(\xmap(\balpha,\bmu), \umap(\balpha,\bmu)).
\end{align}
%
% Along the lines of Figure~\ref{fig:trajman_clean},
%
Figure~\ref{fig:trajman_pronto} provides a graphical interpretation of
PRONTO.
% provides a visual representation of the optimal control problem in the space of curves $(\balpha,\bmu)$. 
%
% At each iteration of the algorithm an update direction (in blue) is sought onto the tangent space 
% to the trajectory manifold (in green).
% Then, the updated (infeasible) curve is projected back onto the trajectory manifold (in black) by the
% projection operator.
%
Specifically, PRONTO iteratively refines, for all $k>0$,
a tentative solution of problem~\eqref{PR:compact_problem} according to the update
\begin{align}\label{PR:curve_update}
  % \begin{bmatrix}
  %         \balpha^{k+1}\\ 
  %         \bmu^{k+1}
  %       \end{bmatrix}
  \begin{bmatrix}
    \bx^{k+1} \\ \bu^{k+1}
  \end{bmatrix}
  & =
    \PP\bigg(
    \underbrace{
    \begin{bmatrix}
	  \bx^{k}\\ 
	  \bu^{k}
  \end{bmatrix} + \stepsizek \begin{bmatrix}
	  \bdeltax^k\\
	  \bdeltau^k
        \end{bmatrix}
  }_{(\balpha^{k+1},\bmu^{k+1})}
  \bigg),
\end{align}
where $\stepsizek \in (0,1]$ is the step-size, while the update direction
$(\bdeltax^k,\bdeltau^k)\in \R^\dimalpha \times \R^\dimmu$ is obtained by
minimizing a quadratic approximation of the cost in~\eqref{PR:compact_problem}
over the tangent space $\tanTTk$ at the current trajectory $(\bx^{k}, \bu^{k})$.
%
% The approach is visually represented in Figure~\ref{fig:trajman_pronto}.
%

The update direction $(\bdeltax^k,\bdeltau^k)$ is obtained as the minimizer
  of the following problem
\begin{align}
	\label{eq:pronto_descent_unc}
	\min_{(\bdeltax, \bdeltau)\in \tanTTk}
	&\; \nabla \cost(\bx^k, \bu^k)\T
	\begin{bmatrix}
		\deltax \\ \deltau
	\end{bmatrix}
	\\
	&+
	\frac{1}{2}
	\begin{bmatrix}
		\deltax \\ \deltau
	\end{bmatrix}\T
	\quadPRONTO
	\begin{bmatrix}
		\deltax \\ \deltau
	\end{bmatrix},
	\nonumber
\end{align}
%$J(\cdot,\cdot)$ 
%
where $\quadPRONTO$ is a square
matrix.
% and the update direction is sought onto the tangent space to the trajectory
% manifold $\tanTTk$ (in green in Figure~\ref{fig:trajman_pronto}).
%
In the pure Newton version of PRONTO, 
$\quadPRONTO$ is the second order derivative
of the reduced problem~\eqref{PR:compact_problem}, including also second order derivatives
of the projection operator, i.e.,
\begin{align}
	\label{eq:quadpronto}
	\hspace{-1.5ex}\quadPRONTO :=& \nabla^2\cost(\bx^k,\bu^k) + \nabla^2 \PP(\balpha^k,\bmu^k) \nabla \cost(\bx^k, \bu^k)
	% \nabla \xmap(\balpha^k,\bmu^k) \nabla^2_{xx} \cost(\bx^k,\bu^k) \nabla \xmap(\balpha^k,\bmu^k)\T \\
	% %
	% &+ \nabla^2 \xmap(\balpha^k,\bmu^k) \nabla_{x} \cost(\bx^k,\bu^k) \\
	% %
	% &+ \nabla \umap(\balpha^k,\bmu^k) \nabla^2_{uu} \cost(\bx^k,\bu^k) \nabla \umap(\balpha^k,\bmu^k)\T\\
	% %
	% &+ \nabla^2 \umap(\balpha^k,\bmu^k) \nabla_{u} \cost(\bx^k,\bu^k) \\
	% %
	% &+ 2 \nabla \xmap(\balpha^k,\bmu^k) \nabla^2_{xu} \cost(\bx^k,\bu^k) \nabla \umap(\balpha^k,\bmu^k)\T 
\end{align}
We refer to~\cite{hauser2002projection} for a detailed discussion.
\begin{remark}
	\label{rem:lower_order_pronto}
	Depending on the choice of $\quadPRONTO$
	some lower-order versions of PRONTO are possible, 
	e.g., setting $\quadPRONTO = I$,
	with $I$ being the identity matrix, 
	we obtain a first-order method.
	Another possibility is to chose $\quadPRONTO$ as the second-order
	derivatives of the cost only.
	\oprocend
\end{remark}

It can be shown that the update direction $(\bdeltax^k,\bdeltau^k)$
% \LS{can be obtained by solving a constrained \textbf{Linear Quadratic (LQ) problem} equivalent to problem~\eqref{eq:pronto_descent_unc}.}
is obtained
solving the \textbf{Linear Quadratic (LQ) problem} 
\begin{align}\hspace{-0.1cm}
	 %(\bdeltax^k,\bdeltau^k) &= 
		%\\
		\min_{\bdeltax,\bdeltau} \: &
		\sum_{t=0}^{T-1} 
		\Bigg( \!\!
		\begin{bmatrix} 
		  a_t^k\\ b_t^k
		\end{bmatrix}\T \!\!
		\begin{bmatrix}
		  \deltax_t \\ \deltau_t
		\end{bmatrix} 
		+\frac{1}{2}
		\begin{bmatrix}
		  \deltax_t \\ \deltau_t 
		\end{bmatrix}\T 
		\!\!
		\begin{bmatrix}
			Q_t^k \!\!\! & \!\! S_t^k \\
			{S_t^{k\top}} \!\!\! & \!\! R_t^k
		\end{bmatrix}
		%H_t^k
		\begin{bmatrix}
			\deltax_t \\ \deltau_t 
		\end{bmatrix}
		\!\!
		\Bigg)
		\nonumber
		\\
		& \hspace{1cm}
		+ {a_T^{k\top}} \deltax_T + \deltax_T\T 
		Q_T^k 
		\deltax_T
		\nonumber
		\\
		\subj \: & \: \deltax_{t+1} = A_t^k \deltax_t + B_t^k \deltau_t, \:\: t \in \horizon{0}{T-1}
		\label{eq:pronto_descent}
		\\
		& \: \deltax_{0} = 0,
		\nonumber
\end{align}
where $Q_t^k \in \R^{\dimx \times \dimu}$, $S_t^k \in \R^{\dimx \times \dimu}$
and $R_t^k \in \R^{\dimu \times \dimu}$ are proper weight matrices, components
of $\quadPRONTO$, while $A_t^k, B_t^k, a_t^k, b_t^k$ follow the shorthand
notation in~\eqref{eq:shorthand}.
Notice that in problem~\eqref{eq:pronto_descent}
a quadratic approximation of the cost of the reduced problem
%$J(\cdot,\cdot)$
about the current iterate $(\bx^{k}, \bu^{k})$ is considered 
and the computed update direction $(\bdeltax^k,\bdeltau^k)$
is constrained to the tangent space of current trajectory $(\bx^k,\bu^k)$,
i.e., the set of trajectories satisfying the linearization of the
nonlinear system dynamics $\dynamics(\cdot,\cdot)$ about $(\bx^k,\bu^k)$.

% In the pure Newton version of PRONTO, $Q_t^k, S_t^k, R_t^k$ are second order derivatives 
% of the reduced problem, including also second order derivatives
% of the projection operator, i.e., of the projection maps.
% $J(\cdot, \cdot)$
%
%

% \LS{

% Finally, the updated curve $(\balpha^{k+1}, \bmu^{k+1})$ is projected
% onto the trajectory manifold such that
% %
% \begin{align*}
% 	\begin{bmatrix}
% 		\bx^{k+1} \\ \bu^{k+1}
% 	\end{bmatrix}
% 	= \PP(\balpha^{k+1}, \bmu^{k+1}).
% \end{align*}

% In Figure~\ref{fig:trajman_pronto} 
% each curve $(\balpha,\bmu)$ is projected
% onto the trajectory manifold $\TT$ (in black) by means of the projection operator $\PP$ (the dashed line).

% }

Algorithm~\ref{alg:PRONTO} recaps the procedure described so far.
\begin{algorithm}[htbp]
\begin{algorithmic}[0]
\caption{PRONTO}
\label{alg:PRONTO}
%\REQUIRE 
%trajectory $(\bx^0, \bu^0)$ with $x_0^0 = \xinit$
\FOR{$k = 0, 1, 2 \ldots$}
%\STATE \IN{compute $K_t^k$}
%
%\STATE \blue{// step 1}
\STATE \textbf{Step 1:} compute descent direction $(\bdeltax^k,\bdeltau^k)$
	% \begin{center}
		by solving the \textbf{LQ problem}~\eqref{eq:pronto_descent}
		% equivalent to~\eqref{eq:pronto_descent_unc}
	% \end{center} 
%
% \begin{align}\hspace{-1cm}
% \label{alg:PRONTO:descent_direction}
% \begin{split}
%   (\bz^k,\bv^k) = 
%   \argmin_{\bdeltax,\bdeltau} \: & \:
%   \sum_{t=0}^{T-1} 
%   \begin{bmatrix} 
%     a_t^k\\ b_t^k
%   \end{bmatrix}\T \!\!
%   \begin{bmatrix}
%     \deltax_t \\ \deltau_t
%   \end{bmatrix}  
%   +\frac{1}{2}
%   \begin{bmatrix}
%     \bdeltax \\ \bdeltau 
%   \end{bmatrix}\T \!\!
%   H^k
%   \begin{bmatrix}
%     \bdeltax \\ \bdeltau 
%   \end{bmatrix}  
%   \\
%   \subj \: & \: \deltax_{t+1} = A_t^k \deltax_t + B_t^k \deltau_t, \:\: t \in \horizon{0}{T-1}
%   \\
%   & \: \deltax_{0} = 0
% \end{split}
% \end{align}
% \STATE compute step-size $\stepsizek$
% \STATE set initial condition $x_{0}^{k+1} = \xinit$
\FOR{$t = 0, \dots, T-1$}
%\STATE \blue{// step 2}
\STATE \textbf{Step 2:} update (unfeasible) curve
  \begin{align}
	\label{alg:PRONTO:update_curve}
		\alpha_t^{k+1} & = x_t^k + \stepsizek \deltax_t^k,
		&
		\mu_t^{k+1} & = u_t^k + \stepsizek \deltau_t^k
  \end{align}
  % \begin{align}
	% \label{alg:PRONTO:update_curve}
	% \begin{split}
	% 	\alpha_t^{k+1} & = x_t^k + \stepsizek \deltax_t^k
	% 	\\
	% 	\mu_t^{k+1} & = u_t^k + \stepsizek \deltau_t^k
	% \end{split}
  % \end{align}
  \vspace{-.5cm}
%\STATE \blue{// step 3}
\STATE \textbf{Step 3:} compute new (feasible) trajectory
  \begin{align*}
    u_t^{k+1} & = \mu_t^{k+1} + K_t (\alpha_t^{k+1} - x_t^{k+1} )
    \\
    x_{t+1}^{k+1} & = \dynamics (x_t^{k+1}, u_t^{k+1})
  \end{align*}
\ENDFOR
\ENDFOR
\end{algorithmic}
\end{algorithm}

%%%%%%%%%%%%%%%%%%%%%%%%%%%%%%%%%%%%%%%%%%%%%%%%%%%%%%%%%%%%%%%%%%%%%%%%%%%%%%%
%%%%%%%%%%%%%%%%%%%%%%%%%%%%%%%%%%%%%%%%%%%%%%%%%%%%%%%%%%%%%%%%%%%%%%%%%%%%%%%

\section{\algname/}
\label{sec:emb_feedback_method}
We are ready to present the main contribution of the
note, namely a general set of first-order approaches, called \algname/, for
numerical optimal control.
We start by describing a pure gradient (or steepest) descent implementation
which we call \gradalgname/.

%%%%%%%%%%%%%%%%%%%%%%%%%%%%%%%%%%%%%%%%%%%%%%%%%%%%%%%%%%%%%%%%%%%%%%%%%%%%%%%
\subsection{{Derivation and Convergence of \gradalgname/}}
%
% Before formally giving the algorithm and stating the convergence result of \gradalgname/
% (Algorithm~\ref{alg:CL_sequential_method}), we provide an intuition of its
% derivation.
%
The founding idea of \algname/ is to formulate and solve an unconstrained
optimization problem as done in the strategy shown in
Section~\ref{sec:open_loop_sequential}.  At the same time, we take also
advantage from the beneficial effects of the state feedback of the projection
operator (cf. \eqref{PR:closed_loop_projection} and \eqref{PR:projection_operator})
used in Section~\ref{sec:pronto}.

The proposed procedure is summarized in Algorithm~\ref{alg:CL_sequential_method},
where we use the shorthand notation in~\eqref{eq:shorthand} and we assume that,
for all $k$, the state-input trajectory is initialized at $x_0^k = \xinit$.
\begin{algorithm}[htbp]
	\begin{algorithmic}[0]
		\caption{\gradalgname/}
		\label{alg:CL_sequential_method}
		%\REQUIRE
		%trajectory $(\bx^0, \bu^0)$ with $x_0^0 = \xinit$
		\FOR{$k = 0, 1, 2 \ldots$}
		%\STATE compute $K_t^k$, for all $t \in \horizon{0}{T-1}$
		\STATE set $\lambda_T^k = \nabla \termcost (x_T^k)$
		\FOR{$t = T-1, \dots, 0$}
		%\STATE \blue{// step 1}
		\STATE \textbf{Step 1:}  compute descent direction
		\begin{subequations}
			\label{alg:CL_sequential_method:descent_direction}
			\begin{align}
				%\begin{split}
					\label{alg:CL_sequential_method:descent_direction:lambda}
					\lambda_t^k &=
					( A_t^k - B_t^k K_t )^\top \lambda_{t+1}^k +
					a_t^k - K_t^{\top} b_t^k
					\\
					\label{alg:CL_sequential_method:descent_direction:v}
					\dmu_t^k &= -B_t^{k \top} \lambda_{t+1}^k - b_{t}^k
					\\
					\label{alg:CL_sequential_method:descent_direction:z}
					\dalpha_t^k &= K_t^{\top} \dmu_t^k
				%\end{split}
			\end{align}
		\end{subequations}
		\ENDFOR
		%\IF {termination criterion satisfied}
		%\STATE return
		%\ENDIF
		%\STATE set
		%\begin{align*}
		%\alpha^k_t &= x^k_t + \stepsize z^k_t, \quad t = 0, \dots, T,\\
		%\mu^k_t &= u^k_t + \stepsize v^k_t, \quad t = 0, \dots, T-1.
		%\end{align*}
% 		\STATE compute step-size $\stepsizek$
		%\STATE compute
		%\begin{align*}
		%  \begin{bmatrix}
		%    \balpha^{k+1}
		%    \\
		%    \bmu^{k+1}
		%  \end{bmatrix}
		%  &
		%  =
		%  \begin{bmatrix}
		%    \balpha^k
		%    \\
		%    \bmu^k
		%  \end{bmatrix}
		%  +
		%  \stepsize
		%  \begin{bmatrix}
		%    \bz^k\\
		%    \bv^k
		%  \end{bmatrix}
		%\end{align*}
		%\STATE set initial condition $x_{0}^{k+1} = \xinit$ and
		%\FOR{$t = 0, \dots, T-1$}
		%\STATE compute
		%\begin{align*}
		%    u_{t}^{k+1} & = \mu_t^{k+1} + K_{t}^k (\alpha_t^{k+1} - x_{t}^{k+1})
		%    \\
		%    x_{t+1}^{k+1} & = f(x_t^{k+1}, u_t^{k+1})
		%\end{align*}
		%\ENDFOR
		\FOR{$t = 0, \dots, T-1$}
		%\STATE \blue{// step 2}
		\STATE \textbf{Step 2:}  update (unfeasible) curve
		% \begin{align}
		% 	\label{alg:CL_sequential_method:curve_update}
		% 	\begin{split}
		% 		\alpha_t^{k+1} & = \alpha_t^k + \stepsizek\, \dalpha_t^k
		% 		\\
		% 		\mu_t^{k+1} & = \mu_t^k + \stepsizek\, \dmu_t^k
		% 	\end{split}
		% \end{align}
		\begin{align}
			\label{alg:CL_sequential_method:curve_update}
				\alpha_t^{k+1} & = \alpha_t^k + \stepsizek\, \dalpha_t^k,
				&
				\mu_t^{k+1} & = \mu_t^k + \stepsizek\, \dmu_t^k
		\end{align}
		\vspace{-.5cm}
		%\STATE \blue{// step 3}
		\STATE \textbf{Step 3:}  compute new (feasible) trajectory
		\begin{align}
			\label{alg:CL_sequential_method:update}
			\begin{split}
				u_t^{k+1} & = \mu_t^{k+1} + K_t (\alpha_t^{k+1} - x_t^{k+1} )
			\\
			x_{t+1}^{k+1} & = \dynamics (x_t^{k+1}, u_t^{k+1})
			\end{split}
		\end{align}
		\ENDFOR
		\ENDFOR
	\end{algorithmic}
\end{algorithm}

By comparing the steps~\eqref{alg:SOL:descent_direction} and~\eqref{alg:SOL:trajectory_update} 
% of Algorithm~\ref{alg:SOL} 
with \eqref{alg:CL_sequential_method:descent_direction} and~\eqref{alg:CL_sequential_method:update}
% the ones of
% Algorithm~\ref{alg:CL_sequential_method} 
one can immediately recognize how the
latter (\gradalgname/) is a closed-loop version of the former.
%
% %
% Before formally stating the convergence result of \gradalgname/
% (Algorithm~\ref{alg:CL_sequential_method}), we provide an intuition of its
% derivation.
% %
% \LS{
% The founding idea of \algname/ is to formulate and solve an unconstrained optimization problem
% as done in the strategy shown in Section~\ref{sec:open_loop_sequential}.
% At the same time, we take also advantage from the beneficial effects of the state feedback 
% of the projection operator used in Section~\ref{sec:pronto}.
% } 

% The founding idea of \algname/ is as follows. First, we recast the constrained problem~\eqref{CL:ocp} into
% its unconstrained reduced form .
% Then, we solve the unconstrained optimization problem by means of a gradient method 
% in which the gradient of the cost function is computed exploiting a proper costate dynamics.

Specifically, by embedding the feedback system \eqref{PR:closed_loop_projection}
(defining the projection operator) into the optimal control
problem~\eqref{eq:ocp_original} one obtains the following optimal control
problem
\begin{subequations}
	\vspace{-3ex}
	\label{CL:ocp}
	\begin{align}\label{CL:ocp:cost}
		\min_{\bx, \bu, \balpha, \bmu}
		\:&\:
		\sum_{t=0}^{T-1}
		\stagecost (x_{t},u_{t})
		+
		\termcost (x_{T})
		\\
		\begin{split}
			\label{CL:ocp:dynamics}
			\subj \: & \: x_{t+1} = \dynamics (x_{t}, u_{t})\\
			& \: u_{t} = \mu_{t} + K_{t} (\alpha_{t} - x_{t}),
			{\quad t \in \horizon{0}{T-1}},
			\\
			& \: x_{0} = \alpha_{0} = \xinit.
		\end{split}
	\end{align}
\end{subequations}
%
%
% It is worth noting that although an accessory constraint has been added, 
% the solution of problem~\eqref{CL:ocp} is not suboptimal with respect to the 
% original problem~\eqref{eq:ocp_original}, i.e., the two problems are equivalent.
% %
% See~\cite{hauser2002projection} for a detailed discussion.
%
%
In order to solve problem~\eqref{CL:ocp}, we adopt the approach described in
Section~\ref{sec:open_loop_sequential}.
We recast problem~\eqref{CL:ocp} in its reduced form by expressing both the
state $x_t$ and the input $u_t$ as functions of a state-input curve
$(\balpha, \bmu)$ via two nonlinear maps
\begin{align}
  \label{eq:cl_maps}
		x_t & = \xmap_t (\balpha,\bmu),
			&
		 	u_t &= \umap_t (\balpha,\bmu),
\end{align}
% \begin{align}
%   \label{eq:cl_maps}
%  	\begin{split}
% 		x_t & = \xmap_t (\balpha,\bmu),
% 			\\
% 		 	u_t &= \umap_t (\balpha,\bmu)
% 	\end{split}
% \end{align}
%
for all $t$. We notice that these maps can be seen as the closed-loop
counterparts of $\xmap_t(\bu)$ in Section~\ref{sec:open_loop_sequential}.%
\footnote{We make a slight abuse of notation by using the same symbol
  $\xmap_t$ as in Section~\ref{sec:open_loop_sequential}.
	Again we omit the map dependence on the initial condition.
	}

Therefore, by exploiting~\eqref{eq:cl_maps}, we can obtain a reduced instance of
problem~\eqref{eq:ocp_original} given by
\begin{align}
	\label{CL:unconstrained_problem}
	\begin{split}
		% \min_{\substack{\balpha \in \R^{\dimalpha} \\\bmu \in \R^{\dimmu}}
		% } 
		\min_{\balpha, \bmu} %\:
		&
		% \:
		% \underbracket{
		\sum_{t=0}^{T-1}
    \stagecost (\xmap_{t}(\balpha, \bmu),\umap_{t}(\balpha,\bmu))
    +
    \termcost (\xmap_{T}(\balpha, \bmu))
		% }
	\\
	&=
	\min_{\balpha, \bmu} \:
	% \min_{
	% 	\substack{\balpha \in \R^{\dimalpha} \\\bmu \in \R^{\dimmu}}} 
	J(\balpha, \bmu),
	\end{split}
\end{align}
% %
which is an unconstrained optimization problem in the variables $\balpha\in \R^{\dimalpha}$ and $\bmu\in \R^{\dimmu}$.

We point out that although the stacks of the maps in~\eqref{eq:cl_maps}
correspond to the projection maps in~\eqref{PR:projection_operator}, in the
resolution of problem~\eqref{CL:unconstrained_problem} we do not need to
evaluate their derivatives.

\begin{figure}[]
	\centering
	\includegraphics{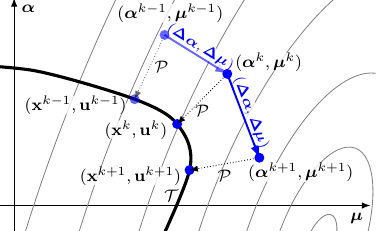}
	\caption{\normalfont 
	Representation of \algname/ approach:
	in gray the level curves of the reduced cost $J(\cdot, \cdot)$, in black
        the trajectory manifold $\TT$, in blue the descent directions.  At each
        iteration $k$, the current curve $(\balpha^k, \bmu^k)$ is updated along
        the (generic) descent direction defined by the gradient of the reduced
        cost $J(\cdot, \cdot)$. The updated curve $(\balpha^{k+1},\bmu^{k+1})$
        is, then, projected onto the trajectory manifold $\TT$ by the projection
        operator $\PP$ (dotted line).}
	%
	%
	% For notational convenience, the update direction is denoted, for all $k$, as $(\dalpha,\dmu)$ instead of $(\dalpha^k,\dmu^k)$.  
	\label{fig:trajman_closedloop}
	\vspace{-4ex}
\end{figure}

Problem~\eqref{CL:unconstrained_problem}, similarly to its open-loop counterpart~\eqref{eq:unconstrained_problem}, 
is an unconstrained optimization problem with nonconvex, twice continuously differentiable cost function $J(\cdot,\cdot)$
(obtained as the composition of $\CC^2$ functions).
Therefore, we apply the gradient method in which the tentative solution $(\balpha^k, \bmu^k)$ is 
iteratively refined as
\begin{align*}
	% \label{eq:CL_alpha_mu}
		\balpha^{k+1}\! =\! \balpha^k 
		\!-\! \stepsizek 
		\nabla_{\balpha}J(\balpha^k,\!\bmu^k),
		\;\;
		\bmu^{k+1}\! = \!\bmu^k 
		\!-\! \stepsizek 
		\nabla_{\bmu} J(\balpha^k,\! \bmu^k)
\end{align*}	
where $k > 0$ is the iteration index while $\stepsizek$ is the step-size.
In parallel with Figure
% ~\ref{fig:trajman_clean} and
~\ref{fig:trajman_pronto},  
a visual representation of this optimization problem is provided 
in Figure~\ref{fig:trajman_closedloop}. 
%
% An intuitive representation of this procedure is given in Figure~\ref{fig:trajman_closedloop}.
We can see that the descent direction is searched in the
entire space of curves $(\balpha, \bmu)$ (rather than on the tangent space to the trajectory
manifold only).  Moreover, the update-direction search is not restricted to any
tangent space.

%  , following a descent direction in the space of the curves.
%
The curve update
% ~\eqref{eq:CL_alpha_mu} 
can be expressed also 
component-wise as
\begin{subequations}
	\label{CL:grad_update}
	\vspace{-3ex}
	\begin{align}
		\alpha_t^{k+1} &= \alpha_t^k 
		- \stepsizek 
		\underbrace{
			\nabla_{\alpha_t}J(\balpha^k, \bmu^k)
			}_{-\dalpha_t^k}
		\label{eq:grad_alpha}
		\\
		\mu_t^{k+1} &= \mu_t^k 
		- \stepsizek 
		\underbrace{
			\nabla_{\mu_t}J(\balpha^k, \bmu^k)
			}_{-\dmu_t^k}
		\label{eq:grad_mu}
		\vspace{-5ex}
	\end{align}
\end{subequations}
for all $t \in \horizon{0}{T-1}$, in which 
each pair $(\dalpha_t^k,\dmu_t^k) \in \R^{\dimx} \times \R^{\dimu}$
represents the descent direction in~\eqref{alg:CL_sequential_method:descent_direction} 
computed by properly adapting the procedure detailed in Section~\ref{sec:open_loop_sequential}.
As it can be seen in Figure~\ref{fig:trajman_closedloop}, each (updated)
state-input curve $(\balpha,\bmu)$ is then projected by the projection operator
$\PP$ onto the trajectory manifold $\TT$ as per \eqref{alg:CL_sequential_method:update}.
% (in black).

% \begin{remark}
% 	We would like to underline that
% 	matrix $K_t$ in~\eqref{alg:CL_sequential_method:descent_direction} and~\eqref{alg:CL_sequential_method:update}, 
% 	for $t \in \horizon{0}{T-1}$, which is designed to 
% 	give stability properties to the projection operator, 
% 	can also be calculated at each iteration $k$,  
% 	%
% 	e.g. by solving the Riccati difference equation associated to a linearization
% 	about current trajectory $(\bx^k, \bmu^k)$.
% 	%
% 	\oprocend
% \end{remark}

Next, we provide the convergence result for Algorithm~\ref{alg:CL_sequential_method}, 
together with an assumption on the step-size.
% Before providing the convergence result for Algorithm~\ref{alg:CL_sequential_method}, 
% let us make an assumption on the step-size.
%
\begin{assumption}
	\label{asm:stepsize}
	Let the step-size $\stepsizek \in \R$, $\stepsizek > 0$ be chosen via
        Armijo backtracking line search.  \oprocend
\end{assumption}
%
%
% The following theorem holds true.
%
\begin{theorem}
	Let Assumptions~\ref{asm:regularity} and \ref{asm:stepsize} hold.
	Let $\{\balpha^k, \bmu^k\}_{k\ge 0}$ be the sequence generated by Algorithm~\ref{alg:CL_sequential_method}. 
	Every limit point $(\balpha^\ast, \bmu^\ast)$ of the sequence $\{\balpha^k, \bmu^k\}_{k\ge 0}$
	satisfies $\nabla J(\balpha^\ast, \bmu^\ast) = 0$.
	Moreover, let $(\bx^\ast, \bu^\ast)$ be the trajectory associated to state-input
	curve $(\balpha^\ast, \bmu^\ast)$ 
	and $\blambda^\ast$ the associated costate trajectory
	%a \IN{limit point NON SONO SICURO} of the sequence of costate 
	%trajectories $\{\blambda^k\}_{k\ge 0}$ 
	generated by Algorithm~\ref{alg:CL_sequential_method}
	in correspondence of $(\balpha^\ast, \bmu^\ast)$.
	Then, $(\bx^\ast, \bu^\ast)$ represents a trajectory satisfying
	the first order necessary conditions for optimality in correspondence of costate trajectory $\blambda^\ast$.
	\oprocend
\label{thm:main_theorem}
\end{theorem}
The proof of Theorem~\ref{thm:main_theorem} can be found in~\cite{sforni2021gopronto}.
\begin{remark}
	Theorem~\ref{thm:main_theorem} can be extended, with suitable assumptions, to different step-size selection rules
	%  other than Armijo backtracking line-search, 
	e.g., constant step-size and diminishing step-size.
 	\oprocend
\end{remark}

% At this point, 
% We observe that since we showed that Algorithm~\ref{alg:CL_sequential_method} corresponds
% to a gradient descent method applied to the reduced form of problem~\eqref{eq:ocp_original},
% it inherits the convergence results of such method.
%

%%%%%%%%%%%%%%%%%%%%%%%%%%%%%%%%%%%%%%%%%%%%%%%%%%%%%%%%%%%%%%%%%%%%%%%%%%%%%%%
%%%%%%%%%%%%%%%%%%%%%%%%%%%%%%%%%%%%%%%%%%%%%%%%%%%%%%%%%%%%%%%%%%%%%%%%%%%%%%%

\subsection{Comparison with Inspiring Methods}
In this subsection we detail the main differences among 
Algorithm~\ref{alg:CL_sequential_method} and the two existing, inspiring algorithms. 
% In this subsection we detail the main differences among 
% Algorithm~\ref{alg:CL_sequential_method} and the two existing, inspiring algorithms. 
%  which gave  to the proposed methodology.

%%%%%%%%%%%%%%%%%%%%%%%%%%%%%%%%%%%%%%%%%%%%%%%%%%%%%%%%%%%%%%%%%%%%%%%%%%%%%%%
%%%%%%%%%%%%%%%%%%%%%%%%%%%%%%%%%%%%%%%%%%%%%%%%%%%%%%%%%%%%%%%%%%%%%%%%%%%%%%%

\subsubsection{Comparison with the gradient method presented in~\cite{bertsekas1999nonlinear}}

\algname/ and the gradient method for optimal control (cf.~Section~\ref{sec:open_loop_sequential}) share the same 
idea of the resolution of the optimal control problem via a gradient method in which
the derivatives are computed through a costate dynamics.
%
% Thanks to the introduction of the projection operator, the consequent
% optimization process takes place along state-input curves rather than input
% sequences only. 
%
However, the introduction of the projection operator implies two fundamental improvements for \algname/.
First, we highlight that \algname/ enjoys numerical stability 
thanks to the different structure of the costate 
dynamics~\eqref{alg:CL_sequential_method:descent_direction:lambda}.
In fact, while 
% in Algorithm~\ref{alg:SOL} 
the dynamical system represented
by~\eqref{alg:SOL:descent_direction_lambda} is governed by the matrix $A_t^k$,
in our algorithm the adjoint system is governed by  the
% the closed-loop 
matrix
$A_t^k - B_t^kK_t$, which for a proper choice of the gain matrices $K_t$
represents a stabilized, time-varying system as the horizon length goes to
infinity.
Moreover, thanks to the projection operator, in \algname/ the input trajectory $\bu^{k+1}$
implements a nonlinear tracking controller of the (updated) state-input 
curve $(\balpha^{k+1}, \bmu^{k+1})$.
Therefore, the trajectory update~\eqref{alg:CL_sequential_method:update} is
performed under a closed-loop strategy rather than in open loop as
in~\eqref{alg:SOL:trajectory_update}, so that dynamical systems subject to
instability issues can be taken into account.

%%%%%%%%%%%%%%%%%%%%%%%%%%%%%%%%%%%%%%%%%%%%%%%%%%%%%%%%%%%%%%%%%%%%%%%%%%%%%%%
%%%%%%%%%%%%%%%%%%%%%%%%%%%%%%%%%%%%%%%%%%%%%%%%%%%%%%%%%%%%%%%%%%%%%%%%%%%%%%%

\subsubsection{Comparison with the PRONTO method presented in~\cite{hauser2002projection}}

\algname/ and PRONTO (cf.~Section~\ref{sec:pronto}) iteratively refine a state-input curve which 
is remapped, using the projection operator, into a state-input trajectory. 
An important difference relies on how these curves are
calculated at each iteration.
In PRONTO, see~\eqref{alg:PRONTO:update_curve}, the next state-input curve is obtained 
by perturbing the current trajectory with a descent direction obtained via 
an LQ problem.
% ~\eqref{eq:pronto_descent}. 
Therefore the direction is sought on the tangent space of the trajectory manifold 
at the current iterate.
%
% Notice that this constrained search space is needed also in the case
% of a first-order implementation of PRONTO (see Remark~\ref{rem:lower_order_pronto} for details).
%
% needs to belong to the tangent space of the trajectory itself, as reflected by the structure of the
% (constrained) LQ problem~\eqref{eq:pronto_descent}.
%
In \algname/, instead, we proceed from curve to curve following the descent
direction defined by the gradient of the reduced 
cost (cf. \eqref{CL:grad_update}) 
%$J(\cdot,\cdot)$
obtained through the adjoint system~\eqref{alg:CL_sequential_method:descent_direction}. 
Since no constraints are imposed on the descent direction and no LQ problems are solved to find the descent direciton, a
lower computation cost is in general required.
%
% Moreover, also sparsity in the problem 
% may be further exploited, see, e.g., \cite{spedicato2018cloud}.

% \medskip

As a final remark, we point out that our algorithmic framework \algname/ 
can be exploited as a globalization technique for Newton's type optimization
methods, see~\cite[Section~1.4]{bertsekas1999nonlinear} for a discussion.
%
%
% can be
% interpreted as a methodological bridge between first-order optimal control
% techniques, as discussed in Section~\ref{sec:open_loop_sequential}, and the
% methods based on stabilizing projection methods, represented by PRONTO.
% %
% % \LS{DOVE METTERE?
% % \begin{remark}
% % 	We would like to underline that
% % 	matrix $K_t$ in~\eqref{alg:CL_sequential_method:descent_direction} and~\eqref{alg:CL_sequential_method:update}, 
% % 	for $t \in \horizon{0}{T-1}$, which is designed to 
% % 	give stability properties to the projection operator, 
% % 	can also be calculated at every iteration $k$,  
% % 	%
% % 	e.g. by solving the Riccati difference equation associated to a linearization
% % 	about the current trajectory $(\bx^k, \bmu^k)$.
% % 	%
% % 	\oprocend
% % \end{remark}
% % }
% %
% In this regard, we point out that our \algname/

%%%%%%%%%%%%%%%%%%%%%%%%%%%%%%%%%%%%%%%%%%%%%%%%%%%%%%%%%%%%%%%%%%%%%%%%%%%%%%%

%%%%%%%%%%%%%%%%%%%%%%%%%%%%%%%%%%%%%%%%%%%%%%%%%%%%%%%%%%%%%%%%%%%%%%%%%%%%%%%
%%%%%%%%%%%%%%%%%%%%%%%%%%%%%%%%%%%%%%%%%%%%%%%%%%%%%%%%%%%%%%%%%%%%%%%%%%%%%%%
%%%%%%%%%%%%%%%%%%%%%%%%%%%%%%%%%%%%%%%%%%%%%%%%%%%%%%%%%%%%%%%%%%%%%%%%%%%%%%%

% \LS{\section{Implementation Details}

% \begin{itemize}
% 	\item Choice of the gain
% 	\item Dealing with constraints
% \end{itemize}}

%%%%%%%%%%%%%%%%%%%%%%%%%%%%%%%%%%%%%%%%%%%%%%%%%%%%%%%%%%%%%%%%%%%%%%%%%%%%%%%

%%%%%%%%%%%%%%%%%%%%%%%%%%%%%%%%%%%%%%%%%%%%%%%%%%%%%%%%%%%%%%%%%%%%%%%%%%%%%%%
%%%%%%%%%%%%%%%%%%%%%%%%%%%%%%%%%%%%%%%%%%%%%%%%%%%%%%%%%%%%%%%%%%%%%%%%%%%%%%%
%%%%%%%%%%%%%%%%%%%%%%%%%%%%%%%%%%%%%%%%%%%%%%%%%%%%%%%%%%%%%%%%%%%%%%%%%%%%%%%

\vspace{-3ex}
\section{Accelerated Versions of \algname/}
\label{sec:accelerated}

In this section we show how 
the framework detailed in
Section~\ref{sec:emb_feedback_method} 
% (thanks to the efficient computation of
% the gradient of the cost function in~\eqref{CL:unconstrained_problem}) 
can be
combined with accelerated gradient-based optimization techniques available in
the literature.
%
%
% This produces accelerated versions of Algorithm~\ref{alg:CL_sequential_method} as described next.
% %
% We would like to underline that these variations
% of \algname/ preserve the dynamic feasibility and numerical stability of the
% original formulation.

%%%%%%%%%%%%%%%%%%%%%%%%%%%%%%%%%%%%%%%%%%%%%%%%%%%%%%%

%%%%%%%%%%%%%%%%%%%%%%%%%%%%%%%%%%%%%%%%%%%%%%%%%%%%%%%%%%%%%%%%%%%%%%%%%%%
%%%%%%%%%%%%%%%%%%%%%%%%%%%%%%%%%%%%%%%%%%%%%%%%%%%%%%%%%%%%%%%%%%%%%%%%%%%

\subsection{\CGalgname/}

% Conjugate Gradient (CG) methods comprise a class of unconstrained optimization
% algorithms.
%which are characterized by low memory requirements and strong local
%and global convergence properties.  
% In the seminal note
% \cite{hestenes1952methods}, the Conjugate Gradient (CG) method is presented as an approach to solve
% symmetric, positive-definite linear systems.  Nevertheless, many strategies were
% developed in the nonlinear system framework, as detailed in
% \cite{hager2006survey}. 
% Details on its implementation are given in
% Appendix~\ref{app:CG}.

The \CGalgname/ optimal control method is obtained by applying the Conjugate Gradient (CG) method to problem
\eqref{CL:unconstrained_problem} (for further details about CG, see~\cite{hestenes1952methods,hager2006survey}).
Let, for all $k>0$ and for all
$t \in \horizon{0}{T-1}$,
\begin{align*}
	% \label{CG:iterate_J}
	\alpha_t^{k+1} &= \alpha_t^k + \stepsizek \tdalpha_t^k, \qquad 
	\mu_t^{k+1} = \mu_t^k +\stepsizek \tdmu_t^k
\end{align*}
where $\stepsizek$ is chosen via Armijo backtracking line search, and
the descent directions $\tdalpha_t^k$ and $\tdmu_t^k$ are obtained according to
the CG algorithm as
%
% \begin{subequations}
	% \label{CG:descent_direction}
\begin{align*}
	\tdalpha_t^k & := 
	%- \nabla_{\alpha_t} J(\balpha^k, \bmu^k) 
	\dalpha_t^k
	+ \CGstep_{\alpha_t}^k \tdalpha_t^{k-1},
	&
	\tdmu_t^k & := 
	%- \nabla_{\mu_t} J(\balpha^k, \bmu^k) 
	\dmu_t^k
	+ \CGstep_{\mu_t}^k \tdmu_t^{k-1}
\end{align*}
% \end{subequations}
%
with 
$\CGstep_{\alpha_t}^k$ and $\CGstep_{\mu_t}^k$ defined as
%
% \begin{subequations}
	\begin{align}
		\label{CG:update_parameter}
		\CGstep_{\alpha_t}^k 
		&:= \tfrac{{\dalpha_t^{k}}\T (\dalpha_t^k - \dalpha_t^{k-1})}
		{\|{\dalpha_t^{k-1} }\|^2},
		%&:= \frac{\nabla_{\alpha_t} J(\balpha^k, \bmu^k)\T (\nabla_{\alpha_t} J(\balpha^k, \bmu^k) - \nabla_{\alpha_t} J(\balpha^{k-1}, \bmu^{k-1}))}
		%{\| \nabla_{\alpha_t} J(\balpha^{k-1}, \bmu^{k-1})\|^2 }
		\hspace{0.5ex}
		\CGstep_{\mu_t}^k
		:=
		\tfrac{{\dmu_t^{k}}\T (\dmu_t^k - \dmu_t^{k-1})}
		{\|{\dmu_t^{k-1} }\|^2}.
		%\frac{\nabla_{\mu_t} J(\balpha^k, \bmu^k)\T (\nabla_{\mu_t} J(\balpha^k, \bmu^k) - \nabla_{\mu_t} J(\balpha^{k-1}, \bmu^{k-1}))}
		%{\| \nabla_{\mu_t} J(\balpha^{k-1}, \bmu^{k-1})\|^2}.
\end{align}
% \end{subequations}
%

Recalling that %the descent direction 
$(\dalpha_t^k, \dmu_t^k)$ (cf.~\eqref{CL:grad_update})
% 
% \begin{subequations}
% \begin{align}
% 	\dalpha_t^k
% 	&= -\nabla_{\alpha_t}J(\balpha^k,\bmu^k)
% 	\\
% 	\dmu_t^k
% 	&= -\nabla_{\mu_t}J(\balpha^k,\bmu^k) 
% \end{align}
% \end{subequations}
can be computed by means of \eqref{alg:CL_sequential_method:descent_direction}, the procedure in
Algorithm~\ref{alg:CG:CL_sequential_method} is obtained.
\begin{algorithm}[htpb]
	\begin{algorithmic}[0]
		\caption{\CGalgname/}
		\label{alg:CG:CL_sequential_method}
		%\REQUIRE
		%trajectory $(\bx^0, \bu^0)$ with $x_0^0 = \xinit$
		\FOR{$k = 0, 1, 2 \ldots$}
		%\STATE compute $K_t^k$, for all $t \in \horizon{0}{T-1}$
		\FOR{$t = T-1, \dots, 0$}
			\STATE \textbf{Step 1:} compute descent direction $\dalpha_t^k, \dmu_t^k$ 
				as in \eqref{alg:CL_sequential_method:descent_direction}
		\ENDFOR

		% \STATE \textbf{Step 1:} compute $\dalpha_t^k, \dmu_t^k$ as 
		% 	in \eqref{alg:CL_sequential_method:descent_direction}
%		\STATE set $\lambda_T^k = \nabla \termcost (x_T^k)$
%		\FOR{$t = T-1, \dots, 0$}
%		\STATE compute
%		\begin{align*}
%			\dmu_t^k &= -B_t^{k \top} \lambda_{t+1}^k - b_{t}^k
%			\\
%			\dalpha_t^k &= K_t^{k \top} \dmu_t^k
%			\\
%			\lambda_t^k &=
%			\Big( A_t^k - B_t^k K_t^k \Big)^\top \lambda_{t+1}^k +
%			a_t^k - K_t^{k \top} b_t^k.
%		\end{align*}
%		\ENDFOR
		\FOR{$t = 0, \dots, T-1$}
		\STATE compute $\CGstep_{\alpha_t}^k, \CGstep_{\mu_t}^k$ as in~\eqref{CG:update_parameter} and the update direction:
		% \begin{align*}
		% 	\umap_{\alpha_t}^k &= \frac{{z_t^{k}}\T (\dalpha_t^k - z_t^{k-1})}
		% 	{\|{z_t^{k-1} }\|^2}
		% 	&
		% 	\umap_{\mu_t}^k &= \frac{{v_t^{k}}\T (\dmu_t^k - v_t^{k-1})}
		% 	{\|{v_t^{k-1} }\|^2}
		% \end{align*}
		% \STATE compute CG update direction:
		\begin{align*}
			\tdalpha_t^k &= 
			\dalpha_t^k
			+ \CGstep_{\alpha_t}^k \tdalpha_t^{k-1}
			&
			\tdmu_t^k &= 
			\dmu_t^k
			+ \CGstep_{\mu_t}^k \tdmu_t^{k-1}
		\end{align*}
		%\ENDFOR
		%\STATE compute step-size $\stepsizek$
		%\FOR{$t = 0, \dots, T-1$}
                \vspace{-3ex}
		\STATE \textbf{Step 2:} update (unfeasible) curve
		% \begin{align*}
		% 	\alpha_t^{k+1} &= \alpha_t^k + \stepsizek \tdalpha_t^k
		% 	\\
		% 	\mu_t^{k+1} &= \mu_t^k + \stepsizek \tdmu_t^k
		% \end{align*}
		\begin{align*}
			\alpha_t^{k+1} &= \alpha_t^k + \stepsizek \tdalpha_t^k,
			&
			\mu_t^{k+1} &= \mu_t^k + \stepsizek \tdmu_t^k
		\end{align*}
                \vspace{-3ex}
		\STATE \textbf{Step 3:} compute new (feasible) trajectory via~\eqref{alg:CL_sequential_method:update}
		% \begin{align*}
		% 	u_t^{k+1} & = \mu_t^{k+1} + K_t^k (\alpha_t^{k+1} - x_t^{k+1} )
		% 	\\
		% 	x_{t+1}^{k+1} & = \dynamics (x_t^{k+1}, u_t^{k+1})
		% \end{align*}
		\ENDFOR
		\ENDFOR
	\end{algorithmic}
\end{algorithm}

As expected, when implemented with the necessary cautions, 
e.g., restarting policies and conjugacy tests,
this method exhibits a faster convergence rate
with respect to its plain gradient counterpart, 
see Section~\ref{sec:simulations} for further details.

%%%%%%%%%%%%%%%%%%%%%%%%%%%%%%%%%%%%%%%%%%%%%%%%%%%%%%%%%%%%%%%%%%%%%%%%%%%

%%%%%%%%%%%%%%%%%%%%%%%%%%%%%%%%%%%%%%%%%%%%%%%%%%%%%%%%%%%%%%%%%%%%%%%%%%%
%%%%%%%%%%%%%%%%%%%%%%%%%%%%%%%%%%%%%%%%%%%%%%%%%%%%%%%%%%%%%%%%%%%%%%%%%%%

\subsection{\HBalgname/}

% The Heavy-Ball method is a two-step procedure for the resolution of
% unconstrained optimization problems.
% % ~\cite{polyak1964some}, 
% %
% It improves the convergence rate with respect to the plain gradient descent.
% %
% Details on its implementation are available in Appendix~\ref{app:HB}.

The \HBalgname/ optimal control method is obtained by applying
the Heavy-ball iteration (cf.~\cite{polyak1964some}) to problem~\eqref{CL:unconstrained_problem}, i.e., for all $k>0$ and for all $t \in \horizon{0}{T-1}$, we have
\begin{align}
	\label{eq:hb_update}
		\begin{split}
			\alpha_t^{k+1} 
			&= \alpha_t^k + \stepsizek\dalpha_t^k
			% \stepsizek \underbrace{\nabla_{\alpha_t}J(\balpha^k,\bmu^k)}_{-\dalpha_t^k} 
			+ \stepsize_{\HB}(\alpha_t^k - \alpha_t^{k-1})
			\\
			\mu_t^{k+1} 
			&= \mu_t^k + \stepsizek\dmu_t^k
			% \stepsizek \underbrace{\nabla_{\mu_t}J(\balpha^k,\bmu^k)}_{-\dmu_t^k}
			+ \stepsize_{\HB}(\mu_t^k - \mu_t^{k-1}),
		\end{split}
\end{align}
where $\stepsizek > 0$ and $\stepsize_{\HB} > 0$ are suitable step-sizes.
% While $\stepsizek$ is usually computed through line-search or minimization rule methods
% the $\stepsize_{\HB}$ is generally constant.
%
The descent directions $\dalpha_t^k, \dmu_t^k$ are computed by means of the
costate equation~\eqref{alg:CL_sequential_method:descent_direction}.
Then, the updated curve obtained via~\eqref{eq:hb_update} is projected into a new (feasible) trajectory via~\eqref{alg:CL_sequential_method:update}.
%
% Algorithm~\ref{alg:HB:CL_sequential_method} recaps the procedure.
% % The resulting procedure is recapped in Algorithm~\ref{alg:HB:CL_sequential_method}.
% %
% \begin{algorithm}[]
% 	\begin{algorithmic}[0]
% 		\caption{\HBalgname/}
% 		\label{alg:HB:CL_sequential_method}
% 		%\REQUIRE
% 		%trajectory $(\bx^0, \bu^0)$ with $x_0^0 = \xinit$
% 		\FOR{$k = 0, 1, 2 \ldots$}
% 		%\STATE compute $K_t^k$, for all $t \in \horizon{0}{T-1}$
% 		\FOR{$t = T-1, \dots, 0$}
% 			\STATE \textbf{Step 1:} compute descent direction $\dalpha_t^k, \dmu_t^k$ 
% 				as in \eqref{alg:CL_sequential_method:descent_direction}
% 		\ENDFOR
% 		\FOR{$t = 0, \dots, T-1$}
% 		\STATE \textbf{Step 2:} update (unfeasible) curve via~\eqref{eq:hb_update}
% 		% \begin{align*}
% 		% 	\alpha_t^{k+1} &= \alpha_t^k + \stepsizek \dalpha_t^k + \stepsize_{\HB}(\alpha_t^k - \alpha_t^{k-1})
% 		% 	\\
% 		% 	\mu_t^{k+1} &= \mu_t^k + \stepsizek \dmu_t^k+ \stepsize_{\HB} (\mu_t^k - \mu_t^{k-1})
% 		% \end{align*}
%                 % \vspace{-3ex}
% 		\STATE \textbf{Step 3:} compute new (feasible) trajectory via~\eqref{alg:CL_sequential_method:update}
% 		% \begin{align*}
% 		% 	u_t^{k+1} & = \mu_t^{k+1} + K_t^k (\alpha_t^{k+1} - x_t^{k+1} )
% 		% 	\\
% 		% 	x_{t+1}^{k+1} & = \dynamics (x_t^{k+1}, u_t^{k+1})
% 		% \end{align*}
% 		\ENDFOR
% 		\ENDFOR
% 	\end{algorithmic}
% \end{algorithm}

We point out that, although a faster convergence rate 
of the Heavy-Ball method (with respect to the plain gradient descent)
is rigorously proved for convex problems only,
the practical implementation of this approach within our methodology
confirmed these expectations (see Section~\ref{sec:simulations}).

%%%%%%%%%%%%%%%%%%%%%%%%%%%%%%%%%%%%%%%%%%%%%%%%%%%%%%%%%%%%%%%%%%%%%%%%%%%

%%%%%%%%%%%%%%%%%%%%%%%%%%%%%%%%%%%%%%%%%%%%%%%%%%%%%%%%%%%%%%%%%%%%%%%%%%%
%%%%%%%%%%%%%%%%%%%%%%%%%%%%%%%%%%%%%%%%%%%%%%%%%%%%%%%%%%%%%%%%%%%%%%%%%%%

\subsection{\nestalgname/}

% Nesterov's accelerated gradient represents an alternative momentum method
% for the resolution of unconstrained optimization problems 
% proposed in~\cite{nesterov1983method}
%
% It is a two step procedure introduced by Nesterov in \cite{nesterov1983method} and
%In general, it is applied to \IN{\bf convex optimization} problems
% 
%
% (details about the generic implementation are available in Appendix~\ref{app:Nesterov})
%
% and has a faster convergence rate than the plain gradient methods.

%
The \nestalgname/ optimal control algorithm is obtained by applying Nesterov's iteration (cf.~\cite{nesterov1983method})
% (cf.~\eqref{app:nest:iter}) 
to problem~\eqref{CL:unconstrained_problem}.
Let, for all $k>0$ and for all $t \in \horizon{0}{T-1}$,
%
% \begin{subequations}
\begin{align}
	\label{Nesterov:update}
	\alpha_t^{k+1} 
	&=
	\tilde{\alpha}_t^k + \stepsizek \tdalpha_t^k,
	% \underbrace{\nabla_{\alpha_t} J(\tbalpha^k,\tbmu^k)}_{-\tdalpha_t^k}
	&
	\mu_t^{k+1} 
	&=
	\tilde{\mu}_t^k + \stepsizek \tdmu_t^k
	% \underbrace{\nabla_{\mu_t} J(\tbalpha^k,\tbmu^k)}_{-\tdmu_t^k},
\end{align}
% \begin{align*}
% 	\alpha_t^{k+1} 
% 	&=
% 	\tilde{\alpha}_t^k + \stepsizek \tdalpha_t^k
% 	% \underbrace{\nabla_{\alpha_t} J(\tbalpha^k,\tbmu^k)}_{-\tdalpha_t^k}
% 	\\
% 	\mu_t^{k+1} 
% 	&=
% 	\tilde{\mu}_t^k + \stepsizek \tdmu_t^k
% 	% \underbrace{\nabla_{\mu_t} J(\tbalpha^k,\tbmu^k)}_{-\tdmu_t^k},
% \end{align*}
% \end{subequations}
% \begin{subequations}
% \begin{align}
% 	\alpha_t^{k+1} 
% 	&=
% 	\tilde{\alpha}_t^k - \stepsizek \underbrace{\nabla_{\alpha_t} J(\tbalpha^k,\tbmu^k)}_{-\tdalpha_t^k}
% 	\\
% 	\mu_t^{k+1} 
% 	&=
% 	\tilde{\mu}_t^k - \stepsizek \underbrace{\nabla_{\mu_t} J(\tbalpha^k,\tbmu^k)}_{-\tdmu_t^k},
% \end{align}
% \end{subequations}
%
where $\stepsizek$ is the step-size, while $\tdalpha_t^k$ and $\tdmu_t^k$
represent the gradient of $J(\cdot,\cdot$) evaluated about the curve $(\tbalpha^k,\tbmu^k)$
% $\tbalpha^k$ and $\tbmu^k$ are the
% stacks of $\tilde{\alpha}_t^k$ and $\tilde{\mu}_t^k$, which are respectively
defined as the stack of
%
% \begin{subequations}
	\begin{align}
		\label{Nesterov:auxiliary_curve}
			\hspace{-1.8ex}\tilde{\alpha}_t^k 
			\!=\!
			\alpha_t^k\! +\! \tfrac{k}{k+3} (\alpha_t^k\! -\! \alpha_t^{k-1}),
			\;\;\tilde{\mu}_t^k 
			\!=\!
			\mu_t^k\! +\! \tfrac{k}{k+3} (\mu_t^k\! -\! \mu_t^{k-1}).
\end{align}
% \end{subequations}
%
% and $\tbalpha^k$ and $\tbmu^k$ the stacks of $\tilde{\alpha}_t^k$ and $\tilde{\mu}_t^k$ respectively.  
%The step-size is generally chosen diminishing \IN{CF REF}. 
%
% Descent directions $\dalpha_t^k, \dmu^t_k$ are calculated, as in Algorithm~\ref{alg:CL_sequential_method}, by means of 
% the adjoint equations~\eqref{alg:CL_sequential_method:descent_direction}. 
%
The procedure is summarized in Algorithm~\ref{alg:Nesterov:CL_sequential_method}.
\begin{algorithm}[]
	\begin{algorithmic}[0]
		\caption{\nestalgname/}
		\label{alg:Nesterov:CL_sequential_method}
		%\REQUIRE
		%trajectory $(\bx^0, \bu^0)$ with $x_0^0 = \xinit$
		\FOR{$k = 0, 1, 2 \ldots$}
		%\STATE compute $K_t^k$, for all $t \in \horizon{0}{T-1}$
		% \STATE \IN{ $\lambda_T^k$  ... }
		\FOR{$t = T-1, \dots, 0$}
		% \FOR{$t = 0, \dots, T-1$}
		\STATE \textbf{Step 1:} compute descent direction
		\begin{align}
			\begin{split}
				\label{Nesterov:adjoint_equations}
				%\label{CL:adjoint_equations:costate}
				\lambda_t^k &=
				\Big( \tilde{A}_t^k - \tilde{B}_t^k K_t \Big)^\top \lambda_{t+1}^k +
				\tilde{a}_t^k - K_t^{\top} \tilde{b}_t^k
				\\
				%\begin{split}
				%\label{CL:adjoint_equations:descent:v}
				\tdmu_t^k &= -\tilde{B}_t^{k \top} \lambda_{t+1}^k - \tilde{b}_{t}^k,
				\quad
				%\label{CL:adjoint_equations:descent:z}
				\tdalpha_t^k = K_t^{\top} \tdmu_t^k.
				%\end{split}
			\end{split}
		\end{align}
		\ENDFOR
		\FOR{$t = 0, \dots, T-1$}
		\STATE compute $\tilde{\alpha}_t^k, \tilde{\mu}_t^k$ as in~\eqref{Nesterov:auxiliary_curve}
		\STATE \textbf{Step 2:} update (unfeasible) curve via~\eqref{Nesterov:update}
		% \begin{align*}
		% 	\alpha_t^{k+1} 
		% 	&=
		% 	\tilde{\alpha}_t^k + \stepsizek\tdalpha_t^k,
		% 	&
		% 	\mu_t^{k+1} 
		% 	&=
		% 	\tilde{\mu}_t^k + \stepsizek\tdmu_t^k
		% \end{align*}
		% \begin{align*}
		% 	\alpha_t^{k+1} 
		% 	&=
		% 	\tilde{\alpha}_t^k + \stepsizek\tdalpha_t^k
		% 	\\
		% 	\mu_t^{k+1} 
		% 	&=
		% 	\tilde{\mu}_t^k + \stepsizek\tdmu_t^k
		% \end{align*}
		\STATE \textbf{Step 3:} compute new (feasible) trajectory via~\eqref{alg:CL_sequential_method:update}
		% \begin{align*}
		% 	u_t^{k+1} & = \mu_t^{k+1} + K_t^k (\alpha_t^{k+1} - x_t^{k+1} )
		% 	\\
		% 	x_{t+1}^{k+1} & = \dynamics (x_t^{k+1}, u_t^{k+1})
		% \end{align*}
		\ENDFOR
		\ENDFOR
	\end{algorithmic}
\end{algorithm}

We point out that the descent direction $(\tdalpha_t^k, \tdmu_t^k)$ in
Algorithm~\ref{alg:Nesterov:CL_sequential_method} is computed at the current
auxiliary curve $(\tbalpha^k, \tbmu^k)$ rather than $(\balpha^k, \bmu^k)$.
% (cf.~\eqref{app:nest:iter:tilde}).
%
% The immediate consequence is that the linearization considered when evaluating
% the adjoint equations is computed about the system trajectory associated to
% $(\tbalpha^k, \tbmu^k)$.
%
In~\eqref{Nesterov:adjoint_equations}, in fact, the
matrices $\tilde{A}_t^k, \tilde{B}_t^k, \tilde{a}_t^k, \tilde{b}_t^k$ 
are defined as
%\begin{subequations}
	\begin{align*}
		\tilde{a}_t^k &:= \nabla_{x_t} \stagecost (\tilde{x}_t^k, \tilde{u}_t^k), \hspace{1.1cm}
		\tilde{b}_t^k := \nabla_{u_t} \stagecost (\tilde{x}_t^k, \tilde{u}_t^k),
		\\
		\tilde{A}_t^k &:= \nabla_{x_t} \dynamics(\tilde{x}_t^k, \tilde{u}_t^k)^\top, \hspace{0.8cm}
		\tilde{B}_t^k := \nabla_{u_t} \dynamics (\tilde{x}_t^k, \tilde{u}_t^k)^\top,
	\end{align*}
%\end{subequations}
with $(\tilde{x}_t^k, \tilde{u}_t^k)= (\xmap_t(\tbalpha^k, \tbmu^k),\umap_t(\tbalpha^k, \tbmu^k))$
for all $k$ and $t$.
%

%
% Finally, we underline also for this method that although being aware its
% convergence results are given only for convex problems, the practical
% implementation of this approach within our framework (see
% Section~\ref{sec:simulations}) presents good results and convergence to
% stationary points of the cost function.
%

%%%%%%%%%%%%%%%%%%%%%%%%%%%%%%%%%%%%%%%%%%%%%%%%%%%%%%%%%%%%%%%%%%%%%%%%%%%%%%%
%%%%%%%%%%%%%%%%%%%%%%%%%%%%%%%%%%%%%%%%%%%%%%%%%%%%%%%%%%%%%%%%%%%%%%%%%%%%%%%

\section{Simulations}
\label{sec:simulations}

\begin{figure}[htpb]
	\centering
	\includegraphics[scale=1]{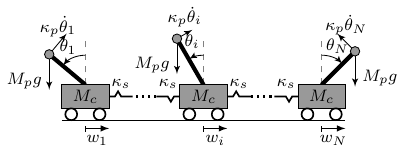}
	\vspace{-2ex}
	\caption{\normalfont Scheme of the train of inverted pendulum-on-cart systems.}
	\label{fig:cart_poles}
\end{figure}

Next, we give explanatory simulations of the algorithms proposed in the previous sections on a large-scale system made by a train of $N$ inverted pendulums on carts
, depicted in Figure~\ref{fig:cart_poles}. 
\LS{The code is implemented in Python on a computer with 2 GHz Quad-Core Intel Core i5 and 16 GB RAM}.
\vspace{-2ex}
\subsubsection{Simulation setup}\label{sec:sim_setup}
For each system $i \in \{1,\ldots, N\}$, the nonlinear continuous-time dynamics is 
\begin{align*}
	M_p l^2 \ddot{\theta}_i + \friction_p\dot{\theta}_i - M_pl\sin(\theta_i)\ddot{w} - M_p l g \sin(\theta_i) 
	& = 0
	\\
	(M_{c} + M_{p})\ddot{w}_i + \friction_{c} \dot{w}_i - \frac{1}{2}M_{p} l \cos(\theta_i)
    \ddot{\theta}_i + \hspace{2cm} &
		\\
		+ \frac{1}{2}M_{p} l \sin(\theta_i) \dot{\theta}_i^2 - \kappa_s w_{i+1} +  \kappa_s w_{i-1} 
		& = u_i,
\end{align*}
where $\theta_i$ is the angle measured from the vertical upward position, 
$w_i$ is the position of the cart and $g =9.81$ m/$\text{s}^2$. The control input is a force $u_i$ applied to the cart.
% Each system is controlled through a force $u_i$ applied to the cart.
%
%
Table~\ref{table:parameters} reports the parameters for all the (identical) carts.
\vspace{-3ex}
\begin{table}[htpb]
	\centering
	\caption{\normalfont Pendulum-on-cart system parameters.\vspace*{-0.2cm}}
	\begin{tabular}{l l l | l l l}
	\multicolumn{3}{c}{Pendulum} & \multicolumn{3}{c}{Cart} \\
	\hline
	Mass & $M_{p}$ & $0.2$ [kg] & Mass & $M_{c}$ & $6.0$ [kg]\\
	Length & $l$ & $1.0$ [m] & Spring Constant & $\kappa_s$ & $0.5$ [$\frac{\text{N}}{\text{m}}$]\\
	Damping & $\friction_{p}$ & $0.01$ [$\frac{\text{Nms}}{\text{rad}}$]
	& Damping & $\friction_{c}$ & $10.0$ [$\frac{\text{Ns}}{\text{m}}$]\\[1.5ex]
	% Gravitational Acceleration & $g$ & $9.81$ [$\frac{\text{m}}{\text{s}^2}$]\\
	\hline
	\end{tabular}
	\label{table:parameters}
	\vspace{-2.5ex}
\end{table}

%
% The state space representation of the dynamics has states 
% $x_i = (\theta_i,\dot{\theta}_i,w_i,\dot{w}_i)\T \in \R^4$ with input $u_i \in \R$ for all $i$.
% %
% Thus, the full state of the system is $x := \col( x_1,\dots,x_{50}) \in \R^{200}$
% and the full input is $u := \col ( u_1, \dots, u_{50}) \in \R^{50}$.
% %
% Then, we use a multiple step Runge-Kutta integrator of order $4$ to obtain a discretized
% version of the plant given by $x_{t+1} = \dynamics(x_t, u_t)$ 
% with sampling period $\delta = 0.05$ seconds.
% The dynamics has states 
% $x_i = (\theta_i,\dot{\theta}_i,w_i,\dot{w}_i)\T \in \R^4$ with input $u_i \in \R$ for all $i$.
%
The discrete-time dynamics, with state $x_{i,t} = (\theta_{i,t},\dot{\theta}_{i,t},w_{i,t},\dot{w}_{i,t})\T \in \R^4$ and input $u_{i,t} \in \R$ for all $i$ is obtained via a Runge-Kutta integrator of order $4$ with sampling period $\delta = 0.05$ seconds.
Thus, the state and the input of the entire system are defined as $x_t := \col( x_{1,t},\dots,x_{N,t}) \in \R^{4N}$
and $u_t := \col ( u_{1,t}, \dots, u_{N,t}) \in \R^{N}$.
% The sensitivities are computed by Algorithmic Differentiation.
%
% For the sake of compactness the discrete time state-space equations are omitted.
We aim at defining the optimal trajectory while tracking a reference curve.
This curve tracking problem has a quadratic cost function %~\eqref{eq:tracking_cost} 
% The tracking problem generally has the following quadratic cost function
\LS{where $\stagecost(x_{t},u_{t}) = \norm{x_t - \xreft}_Q^2 + \norm{u_t - \ureft}_R^2$
and $\termcost(x_T) = \norm{x_T - x_{\text{ref},T}}_{Q_f}^2$}
with symmetric, positive-definite matrices 
$Q := \text{diag}(Q_1, \dots, Q_{N}) \in \R^{4N\times 4N}$ 
and $R := \text{diag}(R_1, \dots, R_{N})$ where, for all $i = 1,\dots, N$,
$Q_i = \text{diag}(100, 1, 0.1, 0.1)$, $R_i = 0.1$.
% \begin{align*}
% 	Q_i &= \begin{bmatrix}
% 		100 & 0 & 0 & 0 \\
% 		0 & 1 & 0 & 0	\\
% 		0 & 0 & 0.1 & 0 \\
% 		0 & 0 & 0 & 0.1
% 	\end{bmatrix}, &
% 	%
% 	R_i &= 10^{-1}.
% \end{align*}
%
The terminal cost matrix $Q_f$ is defined as the solution of the (discrete-time) algebraic Riccati equation evaluated at the linearization of the system about the equilibrium.
We choose a feedback gain $K_t$ in~\eqref{CL:ocp:dynamics} 
solving an LQ problem associated to the linearization of the dynamics 
about the trajectory $(\bx^k, \bu^k)$ available at the current iteration 
with quadratic cost matrices defined as $Q_{\text{reg}}= Q$, $R_{\text{reg}}= I$ and $Q_{f,\text{reg}} = Q_f$.%solution of the algebraic Riccati equation evaluated at the linearization 
%of the system about the equilibrium.
%
\LS{The reference curve is defined, for each cart $i=1,\ldots,N$, as
\begin{align*}
		\theta_\text{i,ref}(t) = \tfrac{\theta_{\text{amp}}^\text{rad} \tanh(t - T/2) (1 - \tanh^2(t - T/2))}{\max\limits_{t \in \horizon{0}{T}}\theta_{\text{amp}}^\text{rad} \tanh(t - T/2) (1 - \tanh^2(t - T/2))}
\end{align*}
where $\theta_{\text{amp}}^\text{rad}$ represent the desired amplitude in radians.
The desired angular velocity is determined by differentiating the smooth curve $\theta_\text{i,ref}(t)$. The other reference signals are zero.}
% carts reference positions, velocities and inputs are zero.}

%%%%%%%%%%%%%%%%%%%%%%%%%%%%%%%%%%%%%%%%%%%%%%%%%%%%%%%%%%%%%%%%%%%%%%%%%%%%%%%%

%%%%%%%%%%%%%%%%%%%%%%%%%%%%%%%%%%%%%%%%%%%%%%%%%%%%%%%%%%%%%%%%%%%%%%%%%%%
%%%%%%%%%%%%%%%%%%%%%%%%%%%%%%%%%%%%%%%%%%%%%%%%%%%%%%%%%%%%%%%%%%%%%%%%%%%
%%%%%%%%%%%%%%%%%%%%%%%%%%%%%%%%%%%%%%%%%%%%%%%%%%%%%%%%%%%%%%%%%%%%%%%%%%%

\subsubsection{Comparison with existing numerical \LS{methods}}

\begin{figure}[t]
	\begin{center}
		\includegraphics{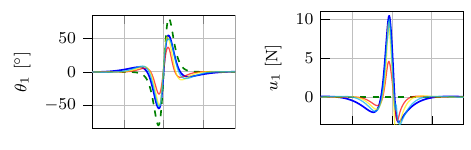}
	\end{center}
	\vspace{-3ex}
	\caption{\normalfont \LS{Optimal angle and input trajectory obtained for the first pendulum-on-cart generated via GoPRONTO with $\theta_{\rm amp} = 80^\circ$.
	In blue the optimal trajectory, in dashed green the reference signals.
	In red, yellow and cyan the trajectories at iteration $k=0,2,4$, respectively.}}
\label{fig:optimal_train}
\vspace{-2ex}
\end{figure}

\LS{In this subsection, we compare \algname/ and the 
SQP-based algorithms for optimal control based on qpOASES, qpDUNES, HPIPM and the first-order OSQP, available in, e.g., \texttt{acados}~\cite{verschueren2022acados}.}
\LS{For all algorithms, we consider the same numerical formulation of the optimal control problem and discretization scheme described above for $N = 3, 50, 100$ systems.}
We aim at performing a swing maneuver between $\pm{\theta_{\rm amp}}$, $\theta_{\rm amp} \in\{60^\circ, 80^\circ\}$ along the smooth reference curve $\theta_\text{i,ref}(t)$ representing the angular reference signal for each $\theta_i$.
\LS{This scenario represents a challenging setting for numerical optimal control algorithms 
due to the large dimension of the decisions variables for the NLP. Moreover, for $\theta_{\rm amp} = 80^\circ$, the problem is even more challenging since it represents an asymptotically ill-posed problem.
We chose this specific setting to provide an insight about the possibilities offered by \algname/,
which  aims at representing a valid alterative in particular scenarios (e.g., quasi ill-posed problems and large-scale systems) where other approaches may face challenge.
%
% The desired angular velocity is determined differentiating the smooth curve for $\theta_i$. 
% Finally, the cart reference positions, velocities and inputs are zero.
%
We remark the fact that \algname/ is implemented as in Algorithm~\ref{alg:CL_sequential_method} without pre-conditioning
nor code optimization.}
\LS{The stepsize is selected via Armijo-line search and the shooting nodes coincide with the discretization time steps.}
\LS{The initial trajectory is $\theta_{i,t} \equiv 0,\dot{\theta}_{i,t} \equiv 0,w_{i,t} \equiv 0,\dot{w}_{i,t} \equiv 0 $, for all $t$ and $i$.}
\LS{The reference signals and the optimal trajectories obtained via \algname/ 
% together with some intermediate ones 
 are shown in Figure~\ref{fig:optimal_train} for the first of $100$ pendulum-on-carts.} 
%  The reference signals are in dashed green.}

% \LS{
	% We compared \algname/ and the SQP-based algorithms for optimal control based on qpOASES, qpDUNES, HPIPM and the first-order OSQP, which are available in, e.g., \texttt{acados}. 
%
% The dynamics is discretized via the same Runge-Kutte scheme used for \algname/. 
% The shooting nodes coincide with the discretization time steps.}
% The initial trajectory is the same both for \algname/ and for the other solvers.}
%
% \LS{\algname/ is compared with the different SQP solvers available in \texttt{acados} which are based on qpOASES, qpDUNES, HPIPM and the first-order OSQP.}
% The \texttt{acados} performances are evaluated adopting 
% We compared \algname with the different QP solvers available,
% i.e., qpOASES, qpDUNES, HPIPM and the first-order OSQP.
%
%
%
\LS{The full condensing versions of qpOASES and HPIPM fail to provide a solution for all $N$ for both references. This can be due to the fact that full condensing approaches eliminate state variables via the unstable dynamics.
%
% Regarding the partial condensing solvers, in which only a few state variables are eliminated, HPIPM successfully solves the problem for $N=3, 50$ for both reference signals.
The partial condensing version of HPIPM, in which only few state variables are eliminated, successfully solves the problem for $N=3, 50$ for both reference signals.
Being HPIPM a second-order solver a faster convergence rate than \algname/ is achieved. 
A comparable convergence rate is achieved only with $\theta_{\rm amp} = 80^\circ$ and $N=50$. This can be due to the significant first-order acceleration required by the ill-posedness of the problem.
% In the latter case, \algname/ achieves a comparable convergence rate, despite HPIPM being a second-order solver.
%
For $N=100$ HPIPM gives segmentation fault. As for partial condensing qpDUNES, it never reaches convergence. % for all $N$.
Finally, we compare \algname/ with the first-order solver OSQP. QSQP has a slower convergence rate for $N=3$ and $\theta_{\rm amp} = 80^\circ$, while it has a faster convergence rate for $\theta_{\rm amp} = 60^\circ$. In the other cases it failes to provide a solution.
Table~\ref{table:acados} summarizes the performances achieved by the solvers that succeded in at least one task.}

\begin{table}[htbp]
	\centering
	\caption{\normalfont Iterations of \gradalgname/, partial condensing HPIPM\\ and OSQP ($\times$ means failure).}
	\setlength\tabcolsep{1.5pt}
	\begin{tabular}{lcccccc}
	\multicolumn{1}{c}{\multirow{2}{*}{$N$}} & \multicolumn{4}{c|}{SQP solvers}                                                            & \multicolumn{2}{c}{\multirow{2}{*}{GoPRONTO}} \\ \cline{2-5}
	\multicolumn{1}{c}{}                     & \multicolumn{2}{c|}{Partial HPIPM}           & \multicolumn{2}{c|}{First-order OSQP}        & \multicolumn{2}{c}{}       \\ \hline
	$\theta_{\rm amp}$                       & $60^\circ$ & \multicolumn{1}{c|}{$80^\circ$} & $60^\circ$ & \multicolumn{1}{c|}{$80^\circ$} & $60^\circ$                 & $80^\circ$       \\ \hline
	$3$                                      & $8$        & $9$                             & $14$       & $50$                            & $30$                       & $34$             \\
	$50$                                     & $14$       & $43$                            & $\times$   & $\times$                        & $35$                       & $39$             \\
	$100$                                    & $\times$   & $\times$                        & $\times$   & $\times$                        & $40$                       & $42$            
	\end{tabular}
	\label{table:acados}
	\vspace{-4ex}
\end{table}

\subsubsection{Comparison with inspiring methods}

\begin{figure}[htpb]
	\centering
	\includegraphics[scale=0.8]{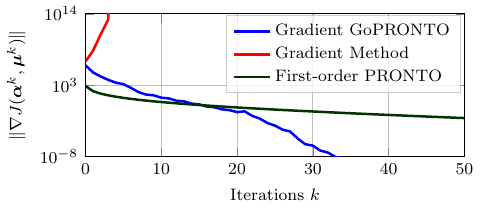}
	\vspace{-1.5ex}
	\caption{\normalfont Evolution of the norm of the gradient $\| \nabla J(\balpha^k,\bmu^k)\|$ 
	in \gradalgname/, the Gradient Method % proposed in~\cite{bertsekas1999nonlinear} 
	and PRONTO.}  % ~\cite{hauser2002projection}
	\label{fig:compare_v_OL_PR}
\end{figure}

Next, \gradalgname/ 
% as detailed in Algorithm~\ref{alg:CL_sequential_method}
is implemented on the previously presented setup with $N = 2$. 
In this case and in the following simulations, the swing maneuver is performed between $+30^\circ$ and $-30^\circ$.
Algorithm~\ref{alg:CL_sequential_method} is compared with the Gradient Method (cf.~Sec.~\ref{sec:open_loop_sequential}) and the first-order version of PRONTO (cf.~Sec.~\ref{sec:pronto}). % As for PRONTO, we implemented the first-order version of the algorithm.
%  where only the second-order derivatives of the cost are considered in~\eqref{eq:quadpronto}.
%
The step-size $\gamma^k$ is selected by Armijo line search rule.
The evolution of the norm of the gradient $\nabla J(\balpha^k,\bmu^k)$, is presented in Figure~\ref{fig:compare_v_OL_PR}.
Notice that the Gradient Method diverges after very few iterations, 
while the first-order version of PRONTO exhibits a slower convergence rate compared with Algorithm~\ref{alg:CL_sequential_method}.
%
% The cost error evolution, represented in Figure~\ref{fig:cost_error}, shows 
% the difference between the cost at iteration $k$ and the asymptotic cost $J^\ast$
% of \algname/.
%
% Notice that the Gradient Method for Optimal Control diverges while both PRONTO and \algname/ converge with almost linear rate. 
% The cost error diminishes across iterations as the optimization proceeds % and new data are gathered 
% with a linear rate, as customary in gradient methods with step-size selected via backtracking approach.
%
% \begin{figure}[htpb]
% 	\centering
% 	\includegraphics[scale=1]{invpend/cost_error.pdf}
% 	\caption{Evolution of the normalized cost error $|\frac{J^k - J^\ast}{J^\ast}|$ in \algname/.} % obtained via DDP}
% 	\label{fig:cost_error}
% \end{figure}
\vspace{-2ex}
\vspace{-2ex}

	\LS{Table~\ref{table:times} compares \algname/ and PRONTO in terms of the computation time required to compute the descent direction and the total computation time with $N = 2,3,5,10$. 
	% These results were obtained on a computer with a 2 GHz Quad-Core Intel Core i5 processor and 16 GB RAM.
	%
	%
	% \LS{The functions required for the implementation of the descent direction calculation routines~\eqref{alg:CL_sequential_method:descent_direction} and~\eqref{eq:pronto_descent}
	% are previously compiled using CasADi.}
	%
	Since PRONTO computes the descent direction by solving an LQ problem, as the state-input dimensions increases, larger computation time per iteration is required with respect to \algname/.
	In \algname/, the total computation time includes the time required to compute the projection gain, which is recomputed only when it looses its stabilizing property, ($2$-$3$ times per simulation).
	While the total computation time is higher with low state dimension, \algname/ shows a faster convergence with high state dimensions.
	Notice that, in PRONTO, one could use the feedback gain obtained by solving the LQ problem at the price of loosing the additional degree of freedom represented by the projection gain.
	It is worth mentioning that GoPRONTO is not a real-time control strategy at the moment, so these results do not provide a comprehensive evaluation of its performance on low-level hardware.}
	\vspace{-2ex}
	\begin{table}[htpb]
		\centering
		\caption{\normalfont Computation times per iteration in [s].\vspace*{-1.5ex}}
			\begin{tabular}{lccccc}
																										& \multicolumn{5}{c}{Computation time per iteration [s]} \\ \cline{2-6} 
			\multicolumn{1}{l|}{Number of carts $N$}      & $2$       & $5$       & $10$      & $50$        & $100$        \\ \hline
			\multicolumn{1}{l|}{PRONTO}                   & $0.12$    & $0.24$    & $0.79$    & $40.86$     & $260.29$     \\
			\multicolumn{1}{l|}{\algname/} & $0.05$    & $0.08$    & $0.18$    & $6.73$      & $33.13$      \\ \hline
																										& \multicolumn{5}{c}{Total computation time [s]}         \\ \hline
			\multicolumn{1}{l|}{PRONTO}                   & $0.84$    & $1.92$    & $7.9$     & $449.6$     & $3123.48$    \\
			\multicolumn{1}{l|}{\algname/} & $1.76$    & $3.36$    & $9.28$    & $328.63$    & $1807.30$   
			\end{tabular}
			\label{table:times}
		\vspace{-2.5ex}
	\end{table}
%%%%%%%%%%%%%%%%%%%%%%%%%%%%%%%%%%%%%%%%%%%%%%%%%%%%%%%%%%%%%%%%%%%%%%%%%%%%%%%%
%%%%%%%%%%%%%%%%%%%%%%%%%%%%%%%%%%%%%%%%%%%%%%%%%%%%%%%%%%%%%%%%%%%%%%%%%%%%%%%%

% \LS{Details about implementation} 
% In Figure~\ref{fig:unc_optimal} the state-input optimal trajectories 
% computed via Algorithm~\ref{alg:CL_sequential_method} are presented.
%
% \begin{figure}[htpb]
% 	\centering
% 	\includegraphics[scale=1]{invpend/trajectory/position.pdf}
% 	% \includegraphics[scale=1]{invpend/trajectory/velocity.pdf}
% 	\includegraphics[scale=1]{invpend/trajectory/input.pdf}
% 	\caption{Optimal trajectory obtained via \algname/. In blue the optimal trajectory, in dashed green the reference signals.}
% 	\label{fig:unc_optimal}
% \end{figure}

%%%%%%%%%%%%%%%%%%%%%%%%%%%%%%%%%%%%%%%%%%%%%%%%%%%%%%%%%%%%%%%%%%%%%%%%%%%%%%%%
%%%%%%%%%%%%%%%%%%%%%%%%%%%%%%%%%%%%%%%%%%%%%%%%%%%%%%%%%%%%%%%%%%%%%%%%%%%%%%%%

\subsubsection{Accelerated versions of \algname/}

In the following, we compare \gradalgname/ with its enhancements.
% : Algorithm~\ref{alg:CG:CL_sequential_method} (\CGalgname/),Algorithm~\ref{alg:HB:CL_sequential_method} (\HBalgname/)
% and Algorithm~\ref{alg:Nesterov:CL_sequential_method} (\nestalgname/). 
%

%%%%%%%%%%%%%%%%%%%%%%%%%%%%%%%%%%%%%%%%%%%%%%%%%%%%%%%%%%%%%%%%%%%%%%%%%%%%%%%%
%%%%%%%%%%%%%%%%%%%%%%%%%%%%%%%%%%%%%%%%%%%%%%%%%%%%%%%%%%%%%%%%%%%%%%%%%%%%%%%%
\subsubsection*{Comparison with \CGalgname/}
\begin{figure}[htpb]
	\vspace{-2ex}
	\centering
	\includegraphics[scale=0.8]{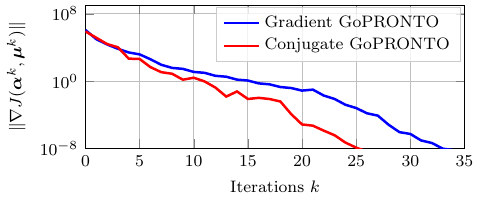}
	\vspace{-1.5ex}
	\caption{\normalfont Evolution of $\|\nabla J(\balpha^k,\bmu^k)\|$ in \gradalgname/ and \CGalgname/ for $N = 2$. $\gamma^k$ is chosen via Armijo line-search.}
	\label{fig:descent_CGvBasic}
\end{figure}

%
% The considered setup is the one detailed previously.
%
Here, $N = 2$ and the step-size $\stepsizek$ is chosen via Armijo line search as required by the Conjugate Gradient method.
Since the CG method is applied to a nonquadratic function, 
we need to deal with the resulting loss of conjugacy.
The implemented method operates in cycles of conjugate direction steps, with the first step of each cycle being a basic gradient step.
We choose to restart the policy when the conjugacy test fails, i.e. as soon as
$|\nabla J(\balpha^{k+1},\bmu^{k+1})\T \nabla J(\balpha^k, \bmu^k)| > 0.7 \| \nabla J(\balpha^k,\bmu^k)\|^2$.
The evolution of the descent direction, i.e., the norm of the gradient $\nabla J(\balpha^k,\bmu^k)$, is presented in Figure~\ref{fig:descent_CGvBasic}.
% In Figure~\ref{fig:descent_CGvBasic} the norm of the gradient $\nabla J(\balpha^k, \bmu^k)$, i.e., the descent direction, is presented.
%
We can see that the descent direction decreases with a faster rate when the CG-enhanced version is adopted.
\vspace{-2ex}

%%%%%%%%%%%%%%%%%%%%%%%%%%%%%%%%%%%%%%%%%%%%%%%%%%%%%%%%%%%%%%%%%%%%%%%%%%%%%%%%

%%%%%%%%%%%%%%%%%%%%%%%%%%%%%%%%%%%%%%%%%%%%%%%%%%%%%%%%%%%%%%%%%%%%%%%%%%%%%%%%
%%%%%%%%%%%%%%%%%%%%%%%%%%%%%%%%%%%%%%%%%%%%%%%%%%%%%%%%%%%%%%%%%%%%%%%%%%%%%%%%
\subsubsection*{Comparison with \HBalgname/ and \nestalgname/}
\begin{figure}[htpb]
	\centering
	\includegraphics[scale=0.8]{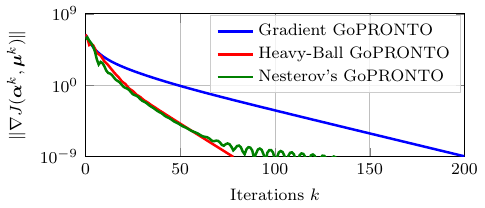}
	\vspace{-1.5ex}
	\caption{\normalfont  Evolution of $\|\nabla J(\balpha^k,\bmu^k)\|$ in \gradalgname/, \nestalgname/ and \HBalgname/ for $N = 50$. The stepsize is constant with $\stepsizek \equiv 10^{-3}$.}
	\label{fig:descent_HBvBasic}
\end{figure}

%
% We consider the setup detailed above.
%
Finally, we consider $N = 50$ resulting in $x_t \in \R^{200}$ and $u_t \in \R^{50}$. The step-size $\stepsizek$ is fixed with $\stepsizek \equiv \stepsize = 10^{-3}$
while the Heavy-ball step $\stepsize_\HB = 0.5$.
%
% In fact, since Nesterov's accelerated gradient is not a descent method, the Armijo backtracking 
% line search is unpracticable. 
%
% We consider in this section also the heavy-ball method for which results about the convergence
% rate are presented under fixed step-size.
%
The evolution of the norm of the gradient $\nabla J(\balpha^k,\bmu^k)$, is presented in Figure~\ref{fig:descent_HBvBasic}.
%
% In Figure~\ref{fig:descent_HBvBasic} the norm of $\nabla J(\balpha^k, \bmu^k)$, i.e., the descent direction, is presented.
%
Notice that the enhanced versions of \algname/ present
a faster convergence rate than its basic implementation.
\vspace{-2ex}

%
% The angular and input optimal trajectories together with some intermediate ones are shown in Figure~\ref{fig:optimal_train} for the first pendulum-on-cart. The reference signals are in dashed green.

%%%%%%%%%%%%%%%%%%%%%%%%%%%%%%%%%%%%%%%%%%%%%%%%%%%%%%%%%%%%%%%%%%%%%%%%%%%%%%%%
%%%%%%%%%%%%%%%%%%%%%%%%%%%%%%%%%%%%%%%%%%%%%%%%%%%%%%%%%%%%%%%%%%%%%%%%%%%%%%%%
%%%%%%%%%%%%%%%%%%%%%%%%%%%%%%%%%%%%%%%%%%%%%%%%%%%%%%%%%%%%%%%%%%%%%%%%%%%%%%%%

\vspace{-2ex}
\subsubsection{Constrained optimal control}
\begin{figure}[htpb]
	\begin{center}
		\includegraphics{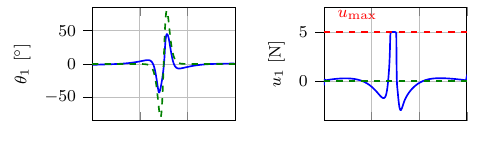}
	\end{center}
	\vspace{-3ex}
	\caption{\normalfont 
\LS{Constrained optimal angle and cart position trajectory for the first pendulum-on-cart.	
	In blue the optimal trajectory, in dashed green the reference trajectory, in dashed red the bound on the control action.}}
\label{fig:optimal_train_const}
\vspace{-2ex}
\end{figure}

\LS{In Figure~\ref{fig:optimal_train_const} we present an example of optimal trajectory where input constraints are enforced. The problem setup is the same as above with $\theta_{\rm amp} = 80^\circ$.
For all carts, i.e., for all $i=1,\ldots, N$, the maximum input is saturated at $u_{\max} = 5$ [N] $g(x_{i,t}) = \| u_{i,t}/u_{\max} \| - 1 \ge 0$, for all $t = \horizon{0}{T}$ 
added to the optimal control problem~\eqref{eq:ocp_original}. This constraint is enforced via the barrier function approach proposed in~\cite{hauser2006barrier}.}
% in Figure~\ref{fig:optimal_train}.
%
% The optimal trajectory for the angular position of some of the 50 carts 
% is represented in Figure~\ref{fig:multiple_theta}.
%
% \begin{figure}[htpb]
% 		\includegraphics{cartpole/trajectory/optimal_traj_2_rev.pdf}
% 		\vspace{-3ex}
% 		\caption{Optimal angle and input trajectory obtained for the first pendulum-on-cart.
% 		In blue the optimal trajectory, in dashed green the reference signals.
% 		In red, yellow and cyan the trajectories at iteration $k=0,2,4$, respectively.}
% 	\label{fig:optimal_train}
% 	\vspace{-2ex}
% \end{figure}
%

%%%%%%%%%%%%%%%%%%%%%%%%%%%%%%%%%%%%%%%%%%%%%%%%%%%%%%%%%%%%%%%%%%%%%%%%%%%

% On this setup, we implemented the plain version of \algname/, i.e, \gradalgname/,
% and its enhanced version \HBalgname/.
% %
% The step-size is constant $\gamma^k \equiv \gamma = 10^{-3}$ while the Heavy-ball step is $\gamma_\HB = 0.5$.
% %
% The evolution of the descent direction, i.e., the norm of the gradient $\nabla J(\balpha^k,\bmu^k)$, is presented in Figure~\ref{fig:descent_train}.
% %
% \begin{figure}[]
% 	\centering
% 	\includegraphics{cartpole/descent_HBvBasic.pdf}
% 	\caption{Evolution of the squared norm of the gradient $\nabla J(\balpha^k,\bmu^k)$ in \gradalgname/ and \HBalgname/ for the $50$ inverted pendulums on carts.}
% 	\label{fig:descent_train}
% \end{figure}
% %
% It is possible to observe that the enhanced version of \algname/ 
% exhibits a faster convergence rate. 

% \begin{figure}[htpb]
% 	\includegraphics{cartpole/multiple_carts/multiple_th.pdf}
% 	\caption{Optimal trajectory for state $\theta_i$ for carts $1,15,30,50$ obtained via \algname/. 
% 	In blue the optimal trajectory, in dashed green the reference signal.}
% 	\label{fig:multiple_theta}
% \end{figure}
% %

%%%%%%%%%%%%%%%%%%%%%%%%%%%%%%%%%%%%%%%%%%%%%%%%%%%%%%%%%%%%%%%%%%%%%%%%%%%
%%%%%%%%%%%%%%%%%%%%%%%%%%%%%%%%%%%%%%%%%%%%%%%%%%%%%%%%%%%%%%%%%%%%%%%%%%%

%%%%%%%%%%%%%%%%%%%%%%%%%%%%%%%%%%%%%%%%%%%%%%%%%%%%%%%%%%%%%%%%%%%%%%%%%%%
%%%%%%%%%%%%%%%%%%%%%%%%%%%%%%%%%%%%%%%%%%%%%%%%%%%%%%%%%%%%%%%%%%%%%%%%%%%
\section{Conclusions}
\label{sec:conclusions}

In this note we proposed \algname/, a novel first-order optimal control methodology that, thanks to the introduction of a nonlinear tracking system, achieves numerical robustness and produces, at each iteration, a feasible trajectory for the system dynamics.
% 
% In the proposed framework, a gradient-based algorithm for optimal control is extended 
% by introducing feedback system (playing the role of a projection operator) in
% the methodology. 
% 
% In this way
% 
Moreover, its simple update rule allowed us to also design several accelerated versions of the plain scheme.
%  as Conjugate gradient, Heavy-ball, and Nesterov's accelerated gradient.
% based on a
% set of adjoint equations.
% % which makes it viable for large-scale dynamical systems. 
% % 
% % In fact, the gradient of the unconstrained problem cost function can be easily calculated
% % by means of a set of adjoint equations.
% % 
% Finally, the gradient-like structure of the proposed framework
%%%%%%%%%%%%%%%%%%%%%%%%%%%%%%%%%%%%%%%%%%%%%%%%%%%%%%%%%%%%%%%%%%%%%%%%%%%%%%%
%%%%%%%%%%%%%%%%%%%%%%%%%%%%%%%%%%%%%%%%%%%%%%%%%%%%%%%%%%%%%%%%%%%%%%%%%%%%%%%

\bibliography{biblio_first_order_closed_loop}
\bibliographystyle{IEEEtran}

%%%%%%%%%%%%%%%%%%%%%%%%%%%%%%%%%%%%%%%%%%%%%%%%%%%%%%%%%%%%%%%%%%%%%%%%%%%%%%%%
%%%%%%%%%%%%%%%%%%%%%%%%%%%%%%%%%%%%%%%%%%%%%%%%%%%%%%%%%%%%%%%%%%%%%%%%%%%%%%%%

\appendices

\section{Proof of Theorem~\ref{thm:main_theorem}}
\label{app:proof}

The proof is arranged in two main parts. In the first part, we prove that any limit point $(\balpha^\ast,\bmu^\ast)$ 
of the sequence $\{\balpha^k,\bmu^k\}_{k\ge0}$ generated by 
Algorithm~\ref{alg:CL_sequential_method} is a stationary point of the unconstrained problem~\eqref{CL:unconstrained_problem},
i.e., it satisfies $\nabla J(\balpha^\ast, \bmu^\ast) = 0$.
Specifically, we show that Algorithm~\ref{alg:CL_sequential_method} represents a 
gradient descent method applied to problem~\eqref{CL:unconstrained_problem}.

Let us prove that the descent direction computed in~\eqref{alg:CL_sequential_method:descent_direction} 
is the gradient of $J(\cdot,\cdot)$ evaluated at the point $(\balpha^k, \bmu^k)$.
%
% \LS{LO ABBIAMO GIA' MESSO PRIMA To this end let us introduce an compact expression for the cost function~\eqref{CL:ocp:cost} 
% %
% \begin{align}
% 	\cost(\bx, \bu) := \sum_{t=0}^{T-1}
% 	\stagecost (x_{t},u_{t})
% 	+
% 	\termcost (x_{T}).
% \label{eq:cost_def}
% \end{align}}
%
To this end let us express the nonlinear dynamics in~\eqref{CL:ocp:dynamics} 
as an implicit equality constraint $\map{\impldynCL}{\R^{\dimalpha } \times \R^{\dimmu } \times \R^{\dimalpha } \times \R^{\dimmu}}{\R^{\dimalpha + \dimmu}}$
defined as
\begin{align}
	\label{eq:impldynCL_def}
	\impldynCL ( \bx,\bu, \balpha,\bmu) :=
	\begin{bmatrix}
		\dynamics (x_0, u_0) - x_{1}\\
		\vdots\\
		\dynamics (x_{T-1}, u_{T-1}) - x_{T}
		%    \\
		%    x_{0} - \xinit
		\\
		\mu_{0} + K_0 (\alpha_0 - x_0) - u_0
		\\
		\vdots\\
		\mu_{T-1} + K_{T-1} (\alpha_{T-1} - x_{T-1}) - u_{T-1}
		%    \\
		%    \alpha_{0} - \xinit
	\end{bmatrix}\!.
\end{align}
We also provide the concise formulation of the cost function~\eqref{eq:ocp_original:cost}
\begin{align}
	\label{eq:cost_def}
	\cost(\bx,\bu) := \sum_{t=0}^{T-1}\stagecost(x_t,u_t) + \termcost(x_T).
\end{align}
Therefore, by means of~\eqref{eq:cost_def}, we can compactly recast problem~\eqref{CL:ocp} as
\begin{align}
	\begin{split}
		\min_{\bx,\bu,\balpha,\bmu} \: & \: \cost(\bx,\bu)
		\\
		\subj \: & \: \impldynCL( \bx,\bu, \balpha,\bmu) = 0.
	\end{split}
	\label{CL:ocp_compact}
\end{align}

Then we can introduce an auxiliary function\footnote{It is evidently the Lagrangian function
of problem~\eqref{CL:ocp_compact}. However, since we do not pursue a Lagrangian approach, 
we prefer not to use such nomenclature.}
associated to problem~\eqref{CL:ocp_compact}, say
$\map{\LL}{\R^\dimalpha \times \R^\dimmu \times R^\dimalpha \times \R^\dimmu \times \R^\dimalpha}{\R}$, 
defined as
\begin{align}
	\LL(\bx,\bu,\balpha,\bmu,\blambda) := \cost(\bx,\bu) + \impldynCL(\bx,\bu,\balpha,\bmu)\T \blambda
\label{CL:lagrangian}
\end{align}
where the (multiplier) vector $\blambda \in \R^{\dimalpha+\dimmu}$ is arranged as
\begin{align*}
	\blambda := \col(\lambda_1,\ldots,\lambda_{T}, \tlambda_{1}, \ldots, \tlambda_{T}  )
\end{align*}
with each $\lambda_t \in \R^\dimx$ and $\tlambda_t \in \R^\dimu$.
By defining $\xmap(\cdot)$ and $\umap(\cdot)$ as the vertical stack 
of the maps $\xmap_t(\cdot)$ and $\umap_t(\cdot)$ (Cf.~\eqref{eq:cl_maps}),
we can see that, by construction, for all $(\balpha,\bmu) \in \R^\dimalpha \times \R^\dimmu$
it holds
\begin{align}\label{CL:h_characterization}
	\impldynCL ( \xmap(\balpha,\bmu),\umap(\balpha,\bmu),\balpha,\bmu ) = 0.
\end{align}
Since $J(\balpha,\bmu) \equiv \cost(\xmap(\balpha,\bmu), \umap(\balpha.\bmu))$
(Cf.~\eqref{CL:unconstrained_problem} and~\eqref{eq:cost_def}), 
the auxiliary function~\eqref{CL:lagrangian} enjoys the following property
\begin{align}
		\label{eq:J_equal_L}
		\LL(\xmap(\balpha,\bmu),\umap(\balpha,\bmu),\balpha,\bmu,\blambda)
		&= J(\balpha,\bmu)
\end{align}
for all $(\balpha,\bmu)$ and \emph{for all} $\blambda\in\R^{\dimalpha+\dimmu}$. 
Therefore, in this formulation we can 
think about $\blambda$ as a parameter or a degree of freedom.
As a consequence of~\eqref{eq:J_equal_L}, it also results
\begin{align}
	\label{eq:gradJ_equal_L}
	\nabla\LL(\xmap(\balpha,\bmu),\umap(\balpha,\bmu),\balpha,\bmu,\blambda)
	&= \nabla J(\balpha,\bmu)
\end{align}
for all $(\balpha, \bmu)$ and, again, \emph{for all} $\blambda$, 
where the gradient of $\LL(\cdot)$ is meant to be calculated only with respect to 
$(\balpha, \bmu)$.
%the variables 
%
%that $\nabla\LL(\cdot) = \nabla J(\cdot)$. 

In the following, we exploit~\eqref{eq:gradJ_equal_L} together with the degree of freedom 
represented by $\blambda$ in order to efficiently compute $\nabla J(\cdot, \cdot)$.
In fact, we can write the two components of the gradient of $J(\cdot,\cdot)$ as
\begin{subequations} \label{eq:grad_J}
\begin{align*}
	\begin{split}
		\nabla_{\balpha} J (\balpha,\bmu)
		& =
		\nabla_{\balpha} \xmap(\balpha,\bmu) \underbracket{\nabla_\bx \LL(\xmap(\balpha,\bmu),\umap(\balpha,\bmu), \balpha,\bmu, \blambda)}
		\\
		&+ \nabla_{\balpha} \umap(\balpha,\bmu) \underbracket{\nabla_\bu \LL(\xmap(\balpha,\bmu),\umap(\balpha,\bmu), \balpha,\bmu, \blambda)}
		\\
		&+ \nabla_{\balpha} \LL(\xmap(\balpha,\bmu),\umap(\balpha,\bmu), \balpha,\bmu, \blambda)
	\end{split}
\end{align*}
and
\begin{align*}
	\begin{split}
		\nabla_{\bmu} J (\balpha,\bmu)
		& =
		\nabla_{\bmu} \xmap(\balpha,\bmu) \underbracket{\nabla_\bx \LL(\xmap(\balpha,\bmu),\umap(\balpha,\bmu), \balpha,\bmu, \blambda)}
		\\
		&+ \nabla_{\bmu} \umap(\balpha,\bmu) \underbracket{\nabla_\bu \LL(\xmap(\balpha,\bmu),\umap(\balpha,\bmu), \balpha,\bmu, \blambda)}
		\\
		&+ \nabla_{\bmu} \LL(\xmap(\balpha,\bmu),\umap(\balpha,\bmu), \balpha,\bmu, \blambda).
	\end{split}
\end{align*}
\end{subequations}
Both these expressions involve the calculation of $\nabla\xmap(\cdot)$ and $\nabla\umap(\cdot)$ 
which may be difficult to compute.
However, since~\eqref{eq:gradJ_equal_L} holds for any $\blambda$, 
we set this degree of freedom to greatly 
simplify the previous formulas.
In fact, the underlined terms $\nabla_{\bx} \LL(\cdot)$ and $\nabla_{\bu}\LL(\cdot)$ have the following peculiar structure
% . 
% By recalling the definition of $\LL(\cdot)$ in~\eqref{CL:lagrangian},
% we can write
%
\begin{subequations}
\label{eq:terms_to_annihilate}
\begin{align}
	\begin{split}
		\nabla_{\bx} \LL(\cdot)
		&= \nabla_\bx \cost ( \xmap(\balpha,\bmu),\umap (\balpha,\bmu)) 
		\\
		& \quad  + \nabla_\bx \impldynCL(\xmap(\balpha,\bmu), \umap (\balpha,\bmu),\balpha,\bmu) \T \blambda 
	\end{split}
\end{align}
and
\begin{align}
	\begin{split}
		\nabla_{\bu} \LL(\cdot)
		&= \nabla_\bu \cost ( \xmap(\balpha,\bmu),\umap (\balpha,\bmu)) 
		\\
		& \quad  + \nabla_\bu \impldynCL(\xmap(\balpha,\bmu), \umap (\balpha,\bmu),\balpha,\bmu) \T \blambda.
	\end{split}
\end{align}
\end{subequations}
Therefore, with a proper choice of $\blambda$ we can annihilate~\eqref{eq:terms_to_annihilate}.
In fact, by choosing $\blambda = \bar{\blambda}$ such that
\begin{align*}
	\nabla_{\bx}\LL(\xmap(\balpha,\bmu),\umap (\balpha,\bmu), \balpha,\bmu, \bar{\blambda})
	&= 0
	\\
	\nabla_{\bu}\LL(\xmap(\balpha,\bmu),\umap (\balpha,\bmu), \balpha,\bmu, \bar{\blambda}) 
	&= 0
\end{align*}
i.e., by setting
\begin{subequations}
	\label{eq:adjoint_system}
	\begin{align}
		\begin{split}
			\label{eq:adjoint_system:x}
			& \nabla_\bx \cost ( \xmap(\balpha,\bmu),\umap (\balpha,\bmu))
			\\
			& \qquad + \nabla_\bx \impldynCL(\xmap(\balpha,\bmu), \umap (\balpha,\bmu),\balpha,\bmu) \T \bar{\blambda} = 0
		\end{split}
		\\
		\begin{split}
			\label{eq:adjoint_system:u}
			& \nabla_\bu \cost ( \xmap(\balpha,\bmu),\umap (\balpha,\bmu))
			\\
			& \qquad + \nabla_\bu \impldynCL(\xmap(\balpha,\bmu), \umap (\balpha,\bmu),\balpha,\bmu) \T \bar{\blambda} = 0,
		\end{split}
	\end{align}
\end{subequations}
both the terms involving $\nabla \xmap(\cdot)$ and $\nabla \umap(\cdot)$ cancel out.
Hence, the gradient components of $J(\cdot,\cdot)$ reduces to
\begin{align*}
	\nabla_{\balpha} J(\balpha,\bmu) 
	&=
	\nabla_{\balpha} \LL(\xmap(\balpha,\bmu),\umap(\balpha,\bmu), \balpha,\bmu, \bar{\blambda})
	\\
	\nabla_{\bmu} J(\balpha,\bmu) 
	&=
	\nabla_{\bmu} \LL(\xmap(\balpha,\bmu),\umap(\balpha,\bmu), \balpha,\bmu, \bar{\blambda}).
\end{align*}

By using again the definition of $\LL(\cdot)$, the latter terms can be written as
\begin{align}
	\label{eq:gradJ_lambda}
	\begin{split}
		\nabla_{\balpha} J (\balpha,\bmu)
	& =
	\nabla_{\balpha} \impldynCL(\xmap(\balpha,\bmu), \umap (\balpha,\bmu),\balpha,\bmu) \T \bar{\blambda}
	\\
	\nabla_{\bmu} J (\balpha,\bmu)
	& =
	\nabla_{\bmu} \impldynCL(\xmap(\balpha,\bmu), \umap (\balpha,\bmu),\balpha,\bmu) \T \bar{\blambda}.
	\end{split}
\end{align}
With this derivation at reach, let us now focus on the $k$-th iteration of Algorithm~\ref{alg:CL_sequential_method}. 
In correspondence of the current state-input curve $(\balpha^k, \bmu^k)$, 
which represents a tentative solution of problem~\eqref{CL:unconstrained_problem}, 
we can compute the vector 
% $\blambda^k \in \R^{\dimalpha+\dimmu}$ arranged as
%
\begin{align*}
	\blambda^k := \col(\lambda_1^k,\ldots,\lambda_{T}^k, \tlambda_{1}^k, \ldots, \tlambda_{T}^k  )
\end{align*}
such that~\eqref{eq:adjoint_system}
holds with $(\balpha, \bmu) = (\balpha^k, \bmu^k)$ and $\bar{\blambda} = \blambda^k$.
Therefore, by recalling the definitions of $\cost(\cdot)$ and $\impldynCL(\cdot)$ in~\eqref{eq:cost_def} and~\eqref{eq:impldynCL_def} 
and since the functions $\dynamics(\cdot), \stagecost(\cdot), \termcost(\cdot)$ are differentiable 
by Assumption~\ref{asm:regularity}, the components $\lambda_t^k$ 
of $\blambda^k$ need to satisfy 
\begin{subequations}
	\begin{align*}
		\nabla \termcost (\xmap_T(\balpha^k,\bmu^k))
		-  \lambda_T^k & = 0
	\end{align*}
and, for all $t \in \horizon{0}{T-1}$,
	\begin{align*}
		a_t^k +
		A_t^{k \top}  \lambda_{t+1}^k
		-\lambda_t^k
		-
		K_t^{\top} \tlambda_{t}^k
		& = 0
	\end{align*}
\end{subequations}
which descends from \eqref{eq:adjoint_system:x}.
As for the components $\tlambda_t^k$ 
of $\blambda^k$, they
needs to be such that for all $t \in \horizon{0}{T-1}$
\begin{align*}
	b_t^k +
	B_t^{k \top}  \lambda_{t+1}^k
	-
	\tlambda_{t}^k
	% \lambda_{t+T}^k
	& = 0
\end{align*}
which comes from \eqref{eq:adjoint_system:u}.
More compactly, a vector $\blambda^k \in \R^{\dimalpha+\dimmu}$
such that for $\bar{\blambda} = \blambda^k$ \eqref{eq:adjoint_system} is satisfied %
for a given $(\balpha^k,\bmu^k)$,
can be obtained by backward simulation of
the adjoint system dynamics
\begin{align}
	\label{eq:adjoint_dynamics}
	\begin{split}
		\lambda_t^k
		& =
		( A_t^k - B_t^k K_t )\T \lambda_{t+1}^k
		+
		a_t^k - K_t^{\top} b_t^k
		\\
		\tlambda_{t}^k
		& =
		b_t^k +
		B_t^{k \top}  \lambda_{t+1}^k
	\end{split}
\end{align}
with terminal condition $\lambda_T^k = \nabla \termcost (\xmap_T(\balpha^k,\bmu^k))$.
With a suitable $\blambda^k$ at hand, we can now compute the gradient of $J(\cdot,\cdot)$
as in~\eqref{eq:gradJ_lambda}, with $(\balpha, \bmu) = (\balpha^k, \bmu^k)$ and
$\bar{\blambda} = \blambda^k$.
Considering a generic time instant $t$ 
and recalling the structure of $\impldynCL(\cdot)$ in~\eqref{eq:impldynCL_def},
we have
\begin{subequations}
\label{eq:gradient_J_final}
\begin{align}
	\label{eq:gradient_alpha}
	\nabla_{\alpha_t} J (\balpha^k,\bmu^k)
	& =
	K_t^{\top} \tlambda_{t}^k  \nonumber
	\\
	& =
	K_t^{\top} \Big( b_{t}^k + B_t^{k \top} \lambda_{t+1}^k \Big)
\end{align}
and
\begin{align}
	\label{eq:gradient_mu}
	\nabla_{\mu_t} J (\balpha^k,\bmu^k)
	& = \tlambda_{t}^k \nonumber
	\\
	& = b_{t}^k + B_t^{k \top} \lambda_{t+1}^k
\end{align}
\end{subequations}
for all $t\in\horizon{0}{T-1}$. 
Comparing \eqref{alg:CL_sequential_method:descent_direction} in Algorithm~\ref{alg:CL_sequential_method} with \eqref{eq:gradient_J_final},
we can see that $\dalpha_t^k, \dmu_t^k$ in \eqref{alg:CL_sequential_method:descent_direction} 
must satisfy
\begin{align*}
	\dalpha_t^k &:= - \nabla_{\alpha_t} J(\balpha^k,\bmu^k)
		%&=- K_t^{\top} \Big( b_{t}^+ B_t^{\top} \lambda_{t+1} \Big)\\
		\\
	\dmu_t^k &:= - \nabla_{\mu_t} J(\balpha^k,\bmu^k).
			%&= - b_{t} - B_t^{\top} \lambda_{t+1}
\end{align*}
Therefore, we proved that Algorithm~\ref{alg:CL_sequential_method} 
tackles problem~\eqref{CL:unconstrained_problem} via a gradient descent method.
%applied to the unconstrained problem~\eqref{CL:unconstrained_problem}.
%\LS{explicitly write gradient descent?}
%
% Since we are using \IN{By ASSUMPTION?} 
In light of Assumption~\ref{asm:stepsize} the step-size $\stepsizek$ in~\eqref{alg:CL_sequential_method:curve_update} is selected according 
to the Armijo rule on the cost function $J(\balpha,\bmu)$.
Therefore, we can conclude that every limit point $(\balpha^*, \bmu^*)$ of $\{\balpha^k,\bmu^k\}_{k\ge0}$ 
is a stationary point of $J(\balpha,\bmu)$, i.e., $\nabla J(\balpha^*,\bmu^*) = 0$.
This completes the first part of the proof.

%%%%%%%%%%%%%%%%%%%%%%%%%%%%%%%%%
%%%%%%%%%%%%%%%%%%%%%%%%%%%%%%%%%
% the point $(\bx^\ast, \bu^\ast)$ we next prove

In the second part, we prove that the state-input trajectory $(\bx^\ast, \bu^\ast) = (\xmap(\balpha^\ast, \bmu^\ast), \umap(\balpha^\ast, \bmu^\ast))$ 
%obtained in correspondence of any limiting state-input curve $(\balpha^\ast, \bmu^\ast)$ is a system trajectory. 
together with the costate vectors $\blambda^\ast \in \R^\dimalpha$ generated by Algorithm~\ref{alg:CL_sequential_method} 
in correspondence of $(\balpha^\ast, \bmu^\ast)$, 
satisfies 
the first order necessary optimality conditions
for the optimal control problem~\eqref{eq:ocp_original}. 
To this end, let us introduce the Hamiltonian function of problem~\eqref{eq:ocp_original} given by
\begin{align*}
	H_t(x_t, u_t, \lambda_{t+1}) := \ell_t(x_t, u_t) + f(x_t, u_t)\T \lambda_{t+1}
\end{align*}
and next we show that
\begin{align*}
	\nabla_{u_t} H_t(x_t^\ast, u_t^\ast, \lambda_{t+1}^\ast) = 0
\end{align*}
and
\begin{align}
	\label{eq:PMP:lambda}
	\lambda_t^\ast = \nabla_{x_t} H_t(x_t^\ast, u_t^\ast, \lambda_{t+1}^\ast)
\end{align}
with terminal condition $\lambda_T^\ast = \nabla \termcost(x_T^\ast)$.

In light of the projection-operator step~\eqref{alg:CL_sequential_method:update}, 
the point $(\bx^\ast,\bu^\ast)$ 
satisfies the dynamics~\eqref{eq:ocp_original:dynamics} by construction, i.e., it is a trajectory.

Let us define the shorthand for the linearization of the cost and 
the dynamics about the trajectory $(\bx^\ast, \bu^\ast)$
\begin{subequations}
	\begin{align}
		a_t^\ast &:= \nabla_{x_t} \stagecost (x_t^\ast, u_t^\ast), \hspace{1.1cm}
		b_t^\ast := \nabla_{u_t} \stagecost (x_t^\ast, u_t^\ast),
		\\
		A_t^\ast &:= \nabla_{x_t} \dynamics(x_t^\ast, u_t^\ast)^\top, \hspace{0.8cm}
		B_t^\ast := \nabla_{u_t} \dynamics (x_t^\ast, u_t^\ast)^\top.
	\end{align}
\end{subequations}
% %
% which represent
%, i.e.
%the system trajectory corresponding to curve $(\balpha^\ast, \bmu^\ast)$.
%
Then, we can define $\blambda^\ast$ as the stack of the costate vectors $\lambda_t^\ast \in \R^\dimx$,
obtained from the adjoint equation~\eqref{alg:CL_sequential_method:descent_direction:lambda} evaluated 
at $(\balpha^\ast, \bmu^\ast)$, i.e., for all $t \in \horizon{T-1}{0}$
\begin{align}
	\label{eq:opt_lambda}
	\lambda_{t}^\ast 
	%=& \; \nabla_{x_t}H_t(x_t^\ast, u_t^\ast, \lambda_{t+1}^\ast)
	%\\
	%=& \; \nabla_{x_t}\stagecost(x_t^\ast, u_t^\ast) - K_t^{\ast\top} \nabla_{u_t}\stagecost(x_t^\ast, u_t^\ast)
	%\\
	%&+ \Big(\nabla_{x_t}\dynamics(x_t^\ast) -\nabla_{u_t}\dynamics(x_t^\ast)K_t^{\ast}\Big)\T \lambda_{t+1}^\ast
	%\\
	=& \; \Big( A_t^\ast - B_t^\ast K_t^\ast \Big)^\top \lambda_{t+1}^\ast + a_t^\ast - K_t^{\top} b_t^\ast
\end{align}
with terminal condition $\lambda_T^\ast = \nabla \termcost(x_T^\ast)$.
Equation~\eqref{eq:opt_lambda} corresponds to the gradient with respect to $x_t$ of the Hamiltonian
evaluated along the trajectory $(\bx^\ast, \bu^\ast)$, i.e., the first order necessary condition for optimality~\eqref{eq:PMP:lambda} holds by construction.

Finally, with $\blambda^\ast$ at hand, we can see that condition
\begin{align*}
	\nabla_{u_t} H_t(x_t^\ast, u_t^\ast, \lambda_{t+1}^\ast) = 0
\end{align*}
can be written as
\begin{align}
	\label{eq:Hamiltonian:nabla_u}
	\nabla_{u_t} H_t(x_t^\ast, u_t^\ast, \lambda_{t+1}^\ast) 
	%&= \nabla_{u_t}\stagecost(x_t^\ast, u_t^\ast) + \nabla_{u_t}\dynamics(x_t^\ast, u_t^\ast)\T \lambda_{t+1}^\ast
	%\\
	&= b_t^\ast + B_t^{\ast\top}\lambda^\ast_{t+1}
\end{align}
which corresponds to $v_t^k$, the gradient of $J(\cdot, \cdot)$ in~\eqref{eq:gradient_mu} evaluated at $(\balpha^\ast, \bmu^\ast)$. 
In light of the first part of the proof, this term is equal to zero.
%, being $(\balpha^\ast, \bmu^\ast)$ a stationary point for $J(\cdot, \cdot)$.
%
Therefore, the first order necessary conditions for optimality are satisfied by
   the trajectory $(\bx^\ast, \bu^\ast)$, thus concluding the proof. 

\end{document}